\documentclass[11pt]{article}%
\usepackage{amsfonts}
\usepackage{amsmath}
\usepackage{amssymb}
\usepackage{amsthm}
\usepackage[a4paper,
            bindingoffset=0.2in,
            left=1in,
            right=1in,
            top=1in,
            bottom=1in,
            footskip=.25in]{geometry}
\usepackage{hyperref}
\usepackage{graphicx}
\usepackage{xspace}
\usepackage{multicol}
\usepackage{contour}%
\usepackage{xcolor}
\usepackage{dsfont}
\usepackage{overpic}
\usepackage{enumerate}
\usepackage{authblk} 
\providecommand{\U}[1]{\protect\rule{.1in}{.1in}}
\theoremstyle{theorem}
\newtheorem{theorem}{Theorem}
\newtheorem{proposition}{Proposition}
\newtheorem{corollary}{Corollary}

\newtheorem*{lemma*}{Lemma}

\theoremstyle{definition}

\newtheorem{remark}{Remark}

\newdimen\dummy
\dummy=\oddsidemargin
\newcommand*{\eg}{e.g.\@\xspace}
\newcommand*{\ie}{i.e.\@\xspace}
\newcommand{\secname}{Section}

\newcommand{\deq}{\mathrel{\mathop:}=}

\contourlength{0.7pt}

\allowdisplaybreaks
\begin{document}

\title{Fractal solutions to the two-dimensional floating body problem}
\author[1]{Yves Martinez-Maure}
\author[2]{David Rochera}
\affil[1]{Sorbonne University and Paris Cit\'{e} University, IMJ-PRG, UMR 7586 of the CNRS, Sophie Germain Building, Box 7012, 75205 Paris Cedex 13, France}
\affil[2]{Dept. of Mathematics, University of Valencia, Av. Vicent Andr\'{e}s Estell\'{e}s 19, E-46100 Burjassot (Valencia), Spain}

\date{}
\setcounter{Maxaffil}{0}
\renewcommand\Affilfont{\itshape\small}

\maketitle

\begin{abstract}
In this paper, the first construction of an entire class of two-dimensional fractal bodies that float in water in equilibrium is explicitly given. Their boundaries are described by planar fractal Zindler curves which are constructed using fractal projective hedgehogs. The areas of these bodies are seen to be well-defined and their expression is explicitly computed. To this end, a result concerning the continuity of the Lebesgue measure with respect to the Hausdorff metric for compact sets is given. Furthermore, the Hausdorff dimension $w$ of these fractal curves is shown to coincide with the Hausdorff dimension of the graph of a Weierstrass function associated to the corresponding fractal projective hedgehog. The fractal Zindler curves are proved to have halving chords of constant length that split both the area and the $w$-dimensional Hausdorff measure of the curve into exact halves.
\end{abstract}

\section{Introduction}

Many mathematical and physical problems related to the equilibrium and stability of a body of uniform density floating in a fluid remain unsolved (or partially unsolved) and have been the subject of numerous research papers over the last decades (see \eg \cite{SchuttWerner90,SchuttWerner94,Stancu06,Stancu09, WernerYe11, HSW19, Gaceta23} and their references therein). One of these problems was posed originally by Ulam (Problem 19 from the Scottish book), who asked whether the sphere is the only (convex) solid of uniform density that floats in water in equilibrium in all orientations.  The density of water is assumed to be equal to 1.

After 80 years, the problem in $\mathbb{R}^n$, for $n\geq 3$, was recently solved by Ryabogin \cite{Ryabogin22annals}, who constructed a strictly convex body of uniform density $\frac{1}{2}$, not a Euclidean ball, that floats in water in equilibrium in every orientation. This means that exactly half of the body's volume is immersed in water, and it remains in equilibrium in every orientation.

The problem in two dimensions is also very interesting. It is known that bodies whose boundaries are Zindler curves are solutions to the floating body problem, providing infinitely many examples different from Euclidean balls. Zindler curves \cite{Zindler21} are curves such that all chords that cut the curve's length into halves have the same length (these are called \emph{halving chords}). Auerbach was the first who used the term \emph{Zindler curve} in \cite{Auerbach38} and studied further properties of these curves. If the halving chords of constant length in a Zindler curve split the corresponding body into exactly two parts, then a remarkable property is that those parts have equal areas. Conversely, if all the chords that divide a smooth plane body into exactly two parts of equal area are of constant length, then the boundary curve of this body is a Zindler curve (whose halving chords are exactly the chords that divide the body area into halves). Such Zindler curves are the boundaries of figures of constant density $\frac{1}{2}$ that float in water in equilibrium in every orientation. There are examples of both convex and non-convex Zindler curves. The general two-dimensional problem for densities $d\neq\frac{1}{2}$ is still open in its full generality. For several densities $d\in\mathopen{}\left]0,\frac{1}{2}\right[\mathclose{}$, there are known examples of these bodies which are not Euclidean balls \cite{Wegner03}. Bracho, Montejano and Oliveros \cite{Carousels04} proved that for perimetral densities $\frac{1}{3}$ and $\frac{1}{4}$, Euclidean balls are the unique solution of the problem, providing thus a positive answer to Ulam's problem.

The objective of this paper is to show that the two-dimensional problem, even for a density $\frac{1}{2}$, exhibits a vast and intricate set of solutions which, remarkably, can be fractal. In fact, we explicitly construct a broad class of fractal Zindler curves that serve as boundaries of figures of constant density $\frac{1}{2}$ that float in water in equilibrium in any position. The construction of these fractal Zindler curves is based on the use of a fractal projective hedgehog.

Classical plane hedgehogs are Minkowski differences of convex bodies in $\mathbb{R}^2$ (see \cite{MM06}). Given a $\mathcal{C}^2$-function $h:\mathbb{S}^1\cong\mathbb{R}/2\pi\mathbb{Z}\to\mathbb{R}$, we can construct the hedgehog with support function $h$ as the envelope $\mathcal{H}_h$ of the family $\bigl(D_{h}(t)\bigr)_{t\in\mathbb{S}^1}$ of support lines
$$x\cos(t)+y\sin(t)=h(t),$$
which yields an explicit parameterization of $\mathcal{H}_h$ in terms of the support function:
\begin{equation}\label{eq:param-hedgehog}
x_h(t)=h(t)\,u(t)+h'(t)\,u'(t),
\end{equation}
where $u(t)=\bigl(\cos(t),\sin(t)\bigr)$. For a detailed introduction to hedgehogs, the interested reader can see \cite{MM-book}.

For any hedgehog $\mathcal{H}_h$, there are exactly two cooriented support lines orthogonal to a given direction $u(t)$, namely $D_h(t)$ and $D_h(t+\pi)$. When these two support lines coincide for every $t\in\mathbb{S}^1$ except for their coorientation (that is, when $h(t)+h(t+\pi)=0$ for all $t\in\mathbb{S}^1$), $\mathcal{H}_h$ is said to be \emph{projective}. It is known that from a projective hedgehog one can construct a family of hedgehogs of constant width and a family of Zindler curves by offsetting a constant distance along the direction given by the coorientation and its orthogonal, respectively (see \eg \cite{MMR22, R22-algebraic, R22-offsets-tire-tracks}).

When the support function is only $\mathcal{C}^1$ the envelope $\mathcal{H}_h$ is still well-defined and parameterized by (\ref{eq:param-hedgehog}), but then $\mathcal{H}_h$ may not correspond to some Minkowski difference of two convex bodies. In fact, it is known that a plane hedgehog that can be regarded as the Minkowski difference of two planar convex bodies is necessarily a rectifiable curve \cite{MM06}. Otherwise, the hedgehog can be fractal as shown in \cite{MM01-fractal}, for instance when $\mathcal{C}^1$-support functions have derivatives that are Weierstrass functions (which are nowhere differentiable).

Based on this idea, in this paper we construct a sequence of smooth Zindler curves such that its (uniform) limit is a fractal Zindler curve. Although this fractal Zindler curve has infinite length, we prove that, provided that the corresponding offsetting is sufficiently large, the curve has no self-intersections, has a finite area and the chords that divide its area into halves have the same length (Theorem~\ref{th:main-theorem-fractal-zindler}). As a consequence, these fractal curves are seen to be the boundaries of figures of constant density $\frac{1}{2}$ that float in water in equilibrium in every orientation (Theorem~\ref{th:fractal-zindler-floating-bodies}).

To ensure the convergence of the sequence of areas of the smooth Zindler curves to the area of the fractal one we need to use a result on the continuity of the Lebesgue measure with respect to the Hausdorff metric for compact sets. Since we have not found this result in the literature and we think it is interesting by itself, we state and prove it here in $\mathbb{R}^d$, with $d$ being a positive integer (Theorem~\ref{th:continuity-lebesgue-measure}).

Since fractal Zindler curves have infinite length, a direct comparison of the two subarcs determined by a chord of constant length cannot be carried out using arc length. Instead, we replace lengths with $w$-dimensional Hausdorff measures, where $w$ is the Hausdorff dimension of the curve, and use this notion to compare the measures of the two resulting fractal subarcs. In this framework, we prove that every halving chord of the fractal Zindler curve divides it into two subarcs of equal $w$-dimensional Hausdorff measure (Corollary~\ref{cor:equality-hauss-measures}). The main result for the Hausdorff measure preservation is stated in Theorem~\ref{th:hauss-measure-preservation}.

The paper is structured as follows. In \secname~\ref{sec:class-curves} we define the class of curves that are used to obtain fractal Zindler curves via a Weierstrass function, that is related to the support function of a fractal projective hedgehog. In \secname~\ref{sec:fractal-zindler}, we prove the first main result. The defined fractal Zindler curves are the boundaries of fractal bodies of constant density $\frac{1}{2}$ that float in water in equilibrium in every orientation. These fractal curves have a well-defined area that is split into halves by chords of constant length, and this area is explicitly computed. In addition, the Hausdorff dimension of the fractal Zindler curve is shown to be equal to the Hausdorff dimension $w$ of the graph of the associated Weierstrass function. In \secname~\ref{sec:hauss-measure-zindler} we prove that the same halving chords also split the $w$-dimensional Hausdorff measure of the fractal Zindler curves into exact halves. Finally, in \secname~\ref{sec:examples} some examples of fractal Zindler curves are shown as well as their relation to fractal hedgehogs of constant width.

\section{A class of fractal curves}\label{sec:class-curves}

First of all, we are going to define the class of curves we will consider throughout the paper. 

Given $n\in\mathbb{N}^{\ast}$, let $\mathcal{H}_{h_n}$ be the projective hedgehog defined by the support function
\begin{equation}\label{eq:supp-func-hedgehogs}
h_n(t)=\sum_{k=1}^n\frac{\sin(\beta^k\,t)}{\alpha^k},
\end{equation}
where $\beta$ is an odd natural number and $\alpha$ is a real number such that $\alpha>\beta$ and $\beta^2>\alpha\,\mathopen{}\left(1+\frac{3\pi}{2}\right)\mathclose{}$. The function $h_n$ is smooth, so a smooth parameterization $x_{h_n}:\mathbb{S}^1\to\mathbb{R}^2$ of $\mathcal{H}_{h_n}$ is naturally given by
$$x_{h_n}(t)=h_n(t)\,u(t)+h_n'(t)\,u'(t),$$
where $u(t)=\bigl(\cos(t),\,\sin(t)\bigr)$. Since $h_n(t)+h_n(t+\pi)=0$ for all $t\in\mathbb{S}^1$, the hedgehog $\mathcal{H}_{h_n}$ is indeed projective.

Notice that the odd natural number $\beta$ must be at least $7$ in order that we could find $\alpha\in\mathbb{R}$ satisfying the conditions $\alpha>\beta$ and $\beta^2>\alpha\,\mathopen{}\left(1+\frac{3\pi}{2}\right)\mathclose{}$. For the chosen values of $\alpha$ and $\beta$, the function
$$h(t)=\lim_{n\to\infty}h_n(t)=\sum_{k=1}^{\infty}\frac{\sin(\beta^k\,t)}{\alpha^k}$$
is of class $\mathcal{C}^1$, but $h'$ is nowhere differentiable (it is a Weierstrass function). The hedgehog $\mathcal{H}_h$ defined by the support function $h$ has a continuous (but nowhere differentiable) parameterization $x_h:\mathbb{S}^1\to\mathbb{R}^2$ given by
$$x_{h}(t)=h(t)\,u(t)+h'(t)\,u'(t).$$
The hedgehog $\mathcal{H}_h$ is a fractal projective hedgehog that has an infinite length \cite{MM01-fractal}. Moreover, the sequence of smooth projective hedgehogs $(\mathcal{H}_{h_n})_{n=1}^{\infty}$ have parameterizations that converge uniformly to that of $\mathcal{H}_h$ on $\mathbb{S}^1$ \cite{MM01-fractal}. 

Let $n\in\mathbb{N}^{\ast}$. Given a projective hedgehog $\mathcal{H}_{h_n}$ and any $r>0$, recall that one can construct a Zindler curve $\mathcal{Z}_{h_n,r}$ by moving each $x_{h_n}(t)$ a constant distance $r$ along the cooriented support line in the direction of $u'(t)$ \cite{MMR22}. The resulting Zindler curve $\mathcal{Z}_{h_n,r}$ associated with $\mathcal{H}_{h_n}$ can be parameterized by $z_{h_n,r}:\mathbb{S}^1\to\mathbb{R}^2$ as
\begin{equation*}
z_{h_n,r}(t)=h_n(t)\,u(t)+\bigl(h_n'(t)+r\bigr)\,u'(t),
\end{equation*}
and its halving chords $\bigl[z_{h_n,r}(t),z_{h_n,r}(t+\pi)\bigr]$ have length $2r$.

\begin{figure}[htb]
\centering
\includegraphics[width=\textwidth]{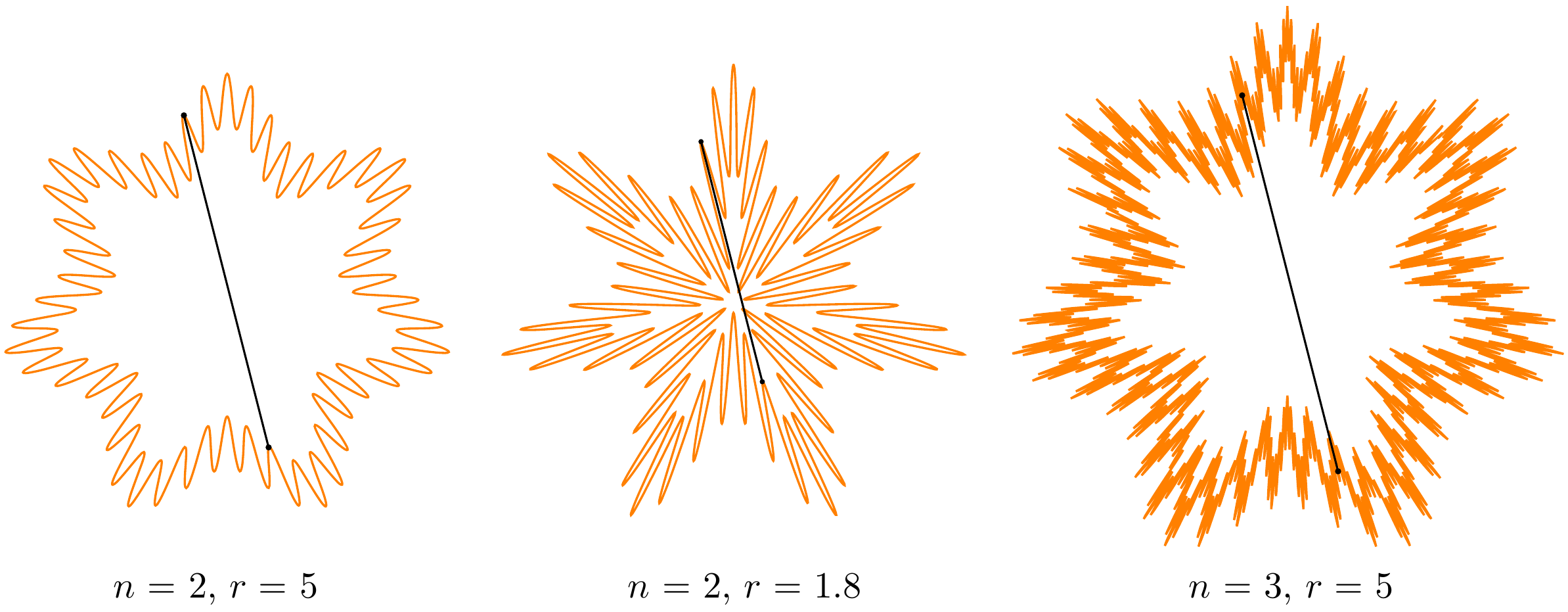}
\caption{Some smooth Zindler curves $\mathcal{Z}_{h_n,r}$ for $\alpha=8$, $\beta=7$ and different values of $n$ and $r$. The halving chord in the middle figure cuts the curve at more than two points.}\label{fig:examples-zindler-r-n}
\end{figure}

Some examples are shown in \figurename~\ref{fig:examples-zindler-r-n}. In particular, notice that the definition of a Zindler curve here refers to curves whose length is split into halves by chords of constant length. This adaptation of the usual definition is necessary because depending on the choice of the parameters $\alpha$, $\beta$, $r$ and $n$, the halving chords could intersect the curve at more than two points (and thus, the set defined by the Zindler curve as its boundary would be split into more than two parts).

\begin{figure}[htb]
\centering
\includegraphics[width=0.45\textwidth]{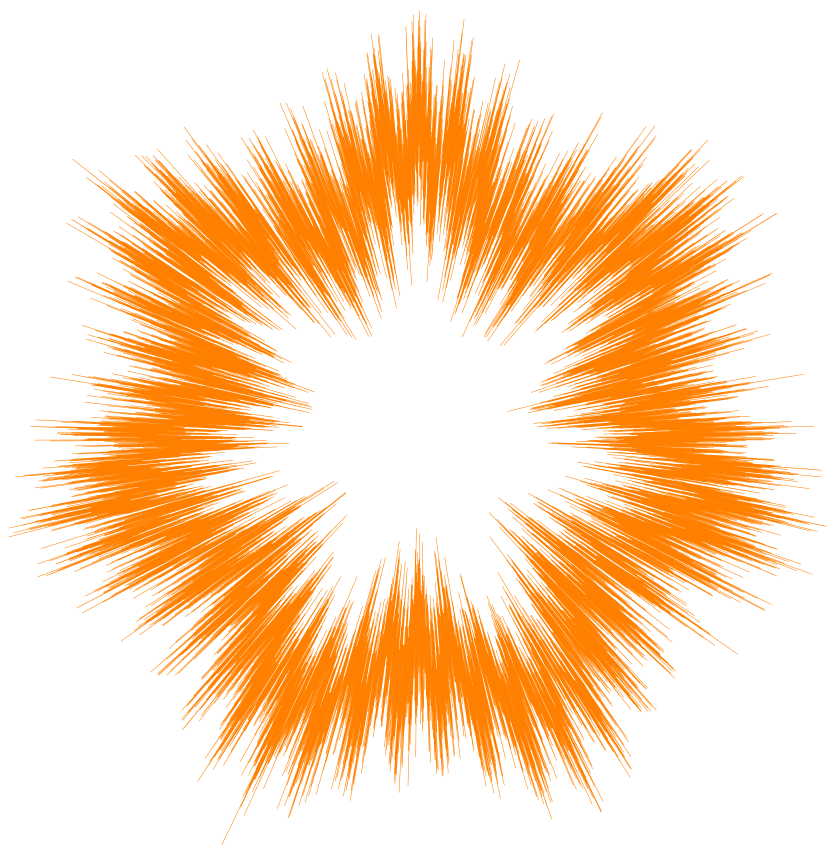}\hspace*{0.05\textwidth}
\includegraphics[width=0.45\textwidth]{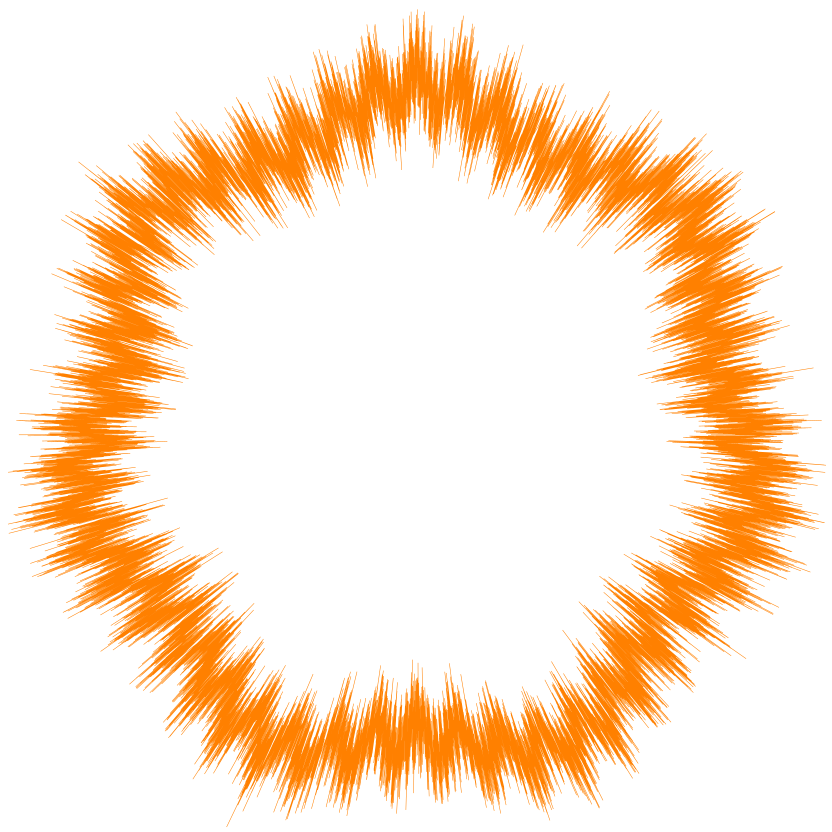}
\caption{An approximation to $\mathcal{Z}_{h,r}$ for $r=7$ (left) and $r=15$ (right), with $\alpha=8$, $\beta=7$ using $n=50$.}\label{fig:fractal-zindler-curves-r}
\end{figure}

Denote by $\mathcal{Z}_{h,r}$ the limit of the sequence of Zindler curves $(\mathcal{Z}_{h_n,r})_{n=1}^{\infty}$ given by the limit of the usual parameterizations of $\mathcal{Z}_{h_n,r}$. That is, $\mathcal{Z}_{h,r}$ is the curve parameterized by $z_{h,r}:\mathbb{S}^1\to\mathbb{R}^2$, where
\begin{equation}\label{eq:usual-param-fractal-zindler-curve}
z_{h,r}(t)=h(t)\,u(t)+\bigl(h'(t)+r\bigr)\,u'(t).
\end{equation}
\figurename~\ref{fig:fractal-zindler-curves-r} shows an approximation to $\mathcal{Z}_{h,r}$ for two values of $r$.

For any $t\in\mathbb{S}^1$, the halving chord $\bigl[z_{h_n,r}(t),z_{h_n,r}(t+\pi)\bigr]$ of $\mathcal{Z}_{h_n,r}$ converges to the chord $\bigl[z_{h,r}(t),z_{h,r}(t+\pi)]$ of $\mathcal{Z}_{h,r}$. We will say that $\bigl[z_{h,r}(t),z_{h,r}(t+\pi)\bigl]$ is the limit of the halving chords with direction $u'(t)$.

The first objective is to show that $\mathcal{Z}_{h,r}$ is a fractal curve (with infinite length) and that it inherits a Zindler property in terms of halving areas. We start with the first point.

\begin{theorem}
The parameterization $z_{h,r}$ of the curve $\mathcal{Z}_{h,r}$ is continuous but nowhere differentiable and $\mathcal{Z}_{h,r}$ has an infinite length.
\end{theorem}

\begin{proof}
The function $h$ is of class $\mathcal{C}^1$ but its derivative $h'$ is nowhere differentiable. As a consequence, the parameterization $z_{h,r}$ of $\mathcal{Z}_{h,r}$ given by (\ref{eq:usual-param-fractal-zindler-curve}) is continuous but nowhere differentiable.

Assume that $\mathcal{Z}_{h,r}$ has a finite length. Then the coordinates of $\mathcal{Z}_{h,r}$ in the canonical basis are necessarily functions of bounded variation. Since functions of bounded variation are almost everywhere differentiable, we find a contradiction with the fact that $z_{h,r}$ is nowhere differentiable. Thus, the length of $\mathcal{Z}_{h,r}$ is infinite.
\end{proof}

Notice that the pointwise convergence defined above of $z_{h_n,r}$ towards $z_{h,r}$ is in fact uniform on $\mathbb{S}^1$. Indeed, given $t\in\mathbb{S}^1$, we have
\begin{align*}
\bigl\|z_{h_n,r}(t)-z_{h,r}(t)\bigr\|&=\Bigl\|\bigl(h_n(t)-h(t)\bigr)\,u(t)+\bigl(h_n'(t)-h'(t)\bigr)\,u'(t)\Bigr\|\\
&\leq \bigl|h_n(t)-h(t)\bigr|+\bigl|h_n'(t)-h'(t)\bigr|.
\end{align*}
Now, the sequences $\bigl(h_n(t)\bigr)_{n=1}^{\infty}$ and $\bigl(h_n'(t)\bigr)_{n=1}^{\infty}$ converge uniformly to $h(t)$ and $h'(t)$ on $\mathbb{S}^1$, respectively, as a consequence of Weierstrass M-test (with our assumptions on $\alpha$ and $\beta$, which are such that $0<1/\alpha<1$ and $0<\beta/\alpha<1$). Therefore, given $\varepsilon>0$, there exists $n_0\in\mathbb{N}$ such that for all $n\geq n_0$ and for all $t\in\mathbb{S}^1$, we have
$$\bigl|h_n(t)-h(t)\bigr|<\frac{\varepsilon}{2}\qquad\text{and}\qquad\bigl|h_n'(t)-h'(t)\bigr|<\frac{\varepsilon}{2}.$$
Thus, for all $n\geq n_0$ and for all $t\in\mathbb{S}^1$, we have
$$\bigl\|z_{h_n,r}(t)-z_{h,r}(t)\bigr\|<\varepsilon,$$
so the convergence of $\mathcal{Z}_{h_n,r}$ towards $\mathcal{Z}_{h,r}$ is uniform on $\mathbb{S}^1$. As a consequence, we also have that the sequence of these curves converges in the Hausdorff metric $d_H$:
$$\lim_{n\to+\infty}d_H\bigl(\mathcal{Z}_{h_n,r},\mathcal{Z}_{h,r}\bigr)=0.$$

\section{Fractal bodies floating in equilibrium in every orientation}\label{sec:fractal-zindler}

We are interested in determining if $\mathcal{Z}_{h,r}$ satisfies a Zindler property in terms of halving areas. As a first step, we compute the area of the smooth Zindler curves $\mathcal{Z}_{h_n,r}$. The next proof uses the computation of some integrals that are gathered in Lemma (Appendix~\ref{sec:appendix}).
\begin{proposition}\label{prop:area-seq-zindler-curves}
Given $\beta\in\mathbb{Z}$ and $\alpha\in\mathbb{R}\setminus\{0,1\}$ such that $\alpha\neq\beta$, the (algebraic) area of the Zindler curve $\mathcal{Z}_{h_n,r}$ is
$$A(\mathcal{Z}_{h_n,r})=\pi\,r^2+\frac{\pi}{2\,\alpha^{2n}}\,\mathopen{}\left(\frac{1-\alpha^{2n}}{1-\alpha^2}-\frac{\beta^2\,(\alpha^{2n}-\beta^{2n})}{\alpha^2-\beta^2}\right)\mathclose{}.$$
\end{proposition}

\begin{proof}
The algebraic area of $z_{h_n,r}(t)=\bigl(z_1(t),\,z_2(t)\bigr)$ can expressed in the form
$$A(\mathcal{Z}_{h_n,r})=\int_0^{2\pi}\frac{1}{2}\,\bigl(z_1(t)\,z_2'(t)-z_2(t)\,z_1'(t)\bigr)\,dt.$$
A straightforward computation shows that
\begin{align*}
A(\mathcal{Z}_{h_n,r})&=\frac{1}{2}\int_0^{2\pi}\bigl(r^2+h_n^2(t)+r\,h_n'(t)+h_n(t)\,h_n''(t)\bigr)\,dt\\
&=\pi\,r^2+\frac{1}{2}\int_0^{2\pi}\bigl(h_n^2(t)+r\,h_n'(t)+h_n(t)\,h_n''(t)\bigr)\,dt.
\end{align*}
Since $h_n(t)=\sum_{k=1}^n\frac{\sin(\beta^k\,t)}{\alpha^k}$, we obtain
\begin{align*}
A(\mathcal{Z}_{h_n,r})=\pi\,r^2&+\frac{1}{2}\int_0^{2\pi}\mathopen{}\left(\sum_{k=1}^n\frac{\sin(\beta^k\,t)}{\alpha^k}\right)^2\mathclose{}dt+\frac{r}{2}\int_0^{2\pi}\mathopen{}\left(\frac{\beta}{\alpha}\right)^k\mathclose{}\cos(\beta^k\,t)\,dt\\
&-\frac{1}{2}\int_0^{2\pi}\mathopen{}\left(\sum_{k=1}^n\frac{\sin(\beta^k\,t)}{\alpha^k}\right)\mathclose{}\,\mathopen{}\left(\sum_{k=1}^n\mathopen{}\left(\frac{\beta^2}{\alpha}\right)^k\sin(\beta^k\,t)\right)\mathclose{}\,dt.
\end{align*}
The value of these three integrals are given in Lemma (Appendix~\ref{sec:appendix}). Thus,
\begin{align*}
A(\mathcal{Z}_{h_n,r})&=\pi\,r^2+\frac{\pi}{2\,\alpha^{2n}}\,\sum_{k=0}^{n-1}\alpha^{2k}-\frac{\beta^2\,\pi}{2\,\alpha^{2n}}\,\sum_{k=0}^{n-1}\alpha^{2n-2k-2}\,\beta^{2k}\\
&=\pi\,r^2+\frac{\pi}{2\,\alpha^{2n}}\,\mathopen{}\left(\sum_{k=0}^{n-1}\alpha^{2k}-\beta^2\,\alpha^{2n-2}\sum_{k=0}^{n-1}\mathopen{}\left(\frac{\beta}{\alpha}\right)^{2k}\mathclose{}\right)\mathclose{}\\
&=\pi\,r^2+\frac{\pi}{2\,\alpha^{2n}}\,\mathopen{}\left(\frac{1-\alpha^{2n}}{1-\alpha^2}-\beta^2\,\alpha^{2n-2}\,\frac{1-\mathopen{}\left(\frac{\beta}{\alpha}\right)^{2n}\mathclose{}}{1-\frac{\beta^2}{\alpha^2}}\right)\mathclose{}\\
&=\pi\,r^2+\frac{\pi}{2\,\alpha^{2n}}\,\mathopen{}\left(\frac{1-\alpha^{2n}}{1-\alpha^2}-\frac{\beta^2\,(\alpha^{2n}-\beta^{2n})}{\alpha^2-\beta^2}\right)\mathclose{},
\end{align*}
as we wanted to show.
\end{proof}

\begin{remark}
If we denote by $(A_n)_{n=1}^{\infty}$ the sequence of areas of the Zindler curves $\mathcal{Z}_{h_n,r}$, then we have
$$A_{n+1}-A_n=\frac{\pi}{2}\,\frac{1-\beta^{2n+2}}{\alpha^{2n+2}},$$
for any $n\in\mathbb{N}^{\ast}$.

Since $\beta>1$ and $\beta<\alpha$, we have $1-\beta^{2n+2}<0<\alpha^{2n+2}$. Thus, $\frac{1-\beta^{2n+2}}{\alpha^{2n+2}}<0$, so that we deduce that  $A_{n+1}-A_n<0$ for all $n\in\mathbb{N}^{\ast}$, which means that the sequence of areas $(A_n)_{n=1}^{\infty}$ is strictly decreasing.
\end{remark}

Once we have computed the area of the curves $\mathcal{Z}_{h_n,r}$, we would like to prove that the area of $\mathcal{Z}_{h,r}$ is the limit of the areas of $\mathcal{Z}_{h_n,r}$ as $n$ goes to infinity. However, there are several issues. The first one is that the curve $\mathcal{Z}_{h,r}$ is a fractal curve and, therefore, it is not trivial to show that it has a well-defined area (an infinite number of self-intersections could happen in a fractal curve). We will prove later in Theorem~\ref{th:main-theorem-fractal-zindler} that indeed the curve $\mathcal{Z}_{h,r}$ defines a Lebesgue measurable set with boundary $\mathcal{Z}_{h,r}$ and has thus a well-defined area. The second issue is that, even in this case, the (plane) Lebesgue measure $\lambda$ is in general not continuous with respect to the Hausdorff metric $d_H$ (it is only upper semi-continuous). This means that, in general, if $K_n$ and $K$ are sets of $\mathbb{R}^2$, $\lim_{n\to+\infty}d_H(K_n,K)=0$ does not imply that $\lambda(K)=\lim_{n\to+\infty}\lambda(K_n)$. However, we will see later that, under our assumptions, it turns out that the Lebesgue measure is continuous with respect to $d_H$. This will be a consequence of the following result.

\begin{theorem}\label{th:continuity-lebesgue-measure}
Let $d$ be a positive integer and let $\lambda$ be the $d$-dimensional Lebesgue measure of $\mathbb{R}^d$. Let $K \subset \mathbb{R}^d$ be a compact set whose boundary $\partial K$ has Hausdorff dimension strictly less than $d$. If a sequence of compact sets $(K_n)_{n=1}^{\infty} \subset \mathbb{R}^d$ converges to $K$ in the Hausdorff metric and the sequence of boundaries $(\partial K_n)_{n=1}^{\infty}$ converges to $\partial K$ also in the Hausdorff metric, then the sequence of Lebesgue measures $\bigl(\lambda(K_n)\bigr)_{n=1}^{\infty}$ converges to $\lambda(K)$.
\end{theorem}

\begin{proof}
Let $\varepsilon>0$. We want to show that there exists $n_0\in\mathbb{N}$ such that for all $n\geq n_0$, $\lambda(K\triangle K_n)<\varepsilon$, where $\triangle$ denotes the symmetric difference of two sets. This implies
$\lim_{n\to+\infty}\lambda(K_n)=\lambda(K)$ using the inequality
$$\bigl|\lambda(K_n)-\lambda(K)\bigr|\leq \lambda(K_n\triangle K),$$
which is easily deduced from basic properties of the Lebesgue measure.

As a first step, let us prove that $\lambda(\partial K)=0$. The Hausdorff dimension of $\partial K$ is defined by
$$\dim_H(\partial K)=\inf\bigl\{\delta\in\mathopen{}\left[0,+\infty\right[\mathclose{}\ :\ \ H^{\delta}(\partial K)=0\bigr\},$$ where $H^{\delta}$ is the $\delta$-dimensional Hausdorff measure. Since $\dim_H(\partial K)<d$, we have $H^d(\partial K)=0$. Now, as $d$ is a positive integer, the $d$-dimensional Hausdorff measure of $\mathbb{R}^d$ is a rescaling of the usual Lebesgue measure $\lambda$ (which is normalized so that the Lebesgue measure of the unit cube $[0,1]^d$ is 1), see \eg \cite{FalconerBook}. In fact, for any Borel subset $X$ of $\mathbb{R}^d$, we have $\lambda(X)=2^{-d}\,v_d\,H^d(X)$, where $v_d$ is the volume of the unit $d$-ball. Therefore, we have $\lambda(\partial K)=0$.

Recall that, given any nonempty set $Y$ of $\mathbb{R}^d$, the distance from a point $x$ in $\mathbb{R}^d$ to $Y$ is defined by
$$d(x,Y)=\inf\bigl\{\|x-y\| :\ y\in Y\bigr\}$$
and that $d(x,Y)=0$ if and only if $x$ is in the closure of $Y$. In addition, given two nonempty sets $X$ and $Y$ of $\mathbb{R}^d$ we have
$$d(X,Y)\deq \inf\bigl\{\|x-y\| :\ x\in X,\ y\in Y\bigr\}=\inf\bigl\{d(x,Y) :\ x\in X\bigr\}.$$

Now, since $\lambda(\partial K)=0$, there exists an open set $U\subset\mathbb{R}^d$ such that $\partial K\subset U$ and $\lambda(U)<\varepsilon$. Note that $d(\partial K,U^c)=a>0$ since $U^c$ is closed and $\partial K$ compact (the map $x\mapsto d(x,U^c)$ is positive, $1$-Lipschitzian, and hence continuous on the compact $\partial K$).

By the assumptions $$\lim_{n\to+\infty}d_H(K_n,K)=0\qquad\text{and}\qquad \lim_{n\to+\infty}d_H(\partial K_n,\partial K)=0,$$
where $d_H$ is the Hausdorff metric, there exists some $n_0\in\mathbb{N}$ such that for all $n\geq n_0$ we have:
$$d_H(K_n,K)<\frac{a}{4}\qquad\text{and}\qquad d_H(\partial K_n,\partial K)<\frac{a}{4}.$$

Recall that the Hausdorff distance between two nonempty comptact subsets $X$ and $Y$ of $\mathbb{R}^d$ is given by
$$d_H(X,Y)=\max\mathopen{}\left(\sup_{x\in X} d(x,Y),\, \sup_{y\in Y} d(y,X)\right)\mathclose{}.$$

Let $n\in\mathbb{N}$, $n\geq n_0$, so that $d_H(K_n,K)<\frac{a}{4}$ and $d_H(\partial K_n,\partial K)<\frac{a}{4}$.

Let $x\in K\setminus K_n$. Since  $x\in K$, we have $d(x,K_n)\leq d_H(K,K_n)<\frac{a}{4}$, and since $x\not\in K_n$, $d(x,K_n)=d(x,\partial K_n)$. It follows that
$$d(x,\partial K)\leq d(x,\partial K_n)+d_H(\partial K_n,\partial K)<\frac{a}{4}+\frac{a}{4}=\frac{a}{2}.$$
Therefore, we have proved that
$$K\setminus K_n\subseteq \mathopen{}\left\{x\in K\ : \ \ d(x,\partial K)<\frac{a}{2}\right\}\mathclose{}.$$
Similarly,
$$K_n\setminus K\subseteq \mathopen{}\left\{x_n\in K_n\ : \ \ d(x_n,\partial K_n)<\frac{a}{2}\right\}\mathclose{}.$$
Of course, $\mathopen{}\left\{x\in K\ : \ \ d(x,\partial K)<\frac{a}{2}\right\}\mathclose{}\subset U$ and if $x_n\in K_n\setminus K$, then
$$d(x_n,\partial K)\leq d(x_n,\partial K_n)+d_H(\partial K_n,\partial K)<\frac{a}{2}+\frac{a}{4} < a,$$
so that
$$\mathopen{}\left\{x_n\in K_n\ : \ \ d(x_n,\partial K_n)<\frac{a}{2}\right\}\mathclose{}\subset U.$$
Therefore,
$$K_n\triangle K = (K_n\setminus K)\cup (K\setminus K_n)\subset U,$$
and we conclude
$$\bigl|\lambda(K_n)-\lambda(K)\bigr|\leq\lambda(K_n\triangle K)\leq \lambda(U)<\varepsilon.\eqno\qedhere$$
\end{proof}

To apply the previous result to our case, we need to show that the Hausdorff dimension of $\mathcal{Z}_{h,r}$ is strictly less than 2. In fact, we are able to prove that if $r$ is sufficiently large, the Hausdorff dimension of $\mathcal{Z}_{h,r}$ is equal to the Hausdorff dimension of the graph of the associated Weierstrass function $h'$, which is known to be $w=3-\frac{\ln(\alpha)}{\ln(\beta)}<2$ (see \cite{Shen18}). This can be found in the following statement, one of the main results of the paper.

\begin{theorem}\label{th:main-theorem-fractal-zindler}
Let $\beta$ be an odd natural number and let $\alpha$ be a real number such that $\alpha>\beta$ and $\beta^2>\alpha\,\mathopen{}\left(1+\frac{3\pi}{2}\right)\mathclose{}$. Given $r>0$, let $\mathcal{Z}_{h,r}$ be the curve parameterized by the uniform limit (\ref{eq:usual-param-fractal-zindler-curve}) of the parameterizations of the smooth Zindler curves $\bigl(\mathcal{Z}_{h_n,r}\bigr)_{n=1}^{\infty}$ associated with the projective hedgehogs $\mathcal{H}_{h_n}$, with $h_n$ as in (\ref{eq:supp-func-hedgehogs}), and with halving chords of length $2r$. If $r$ is sufficiently large, then:
\begin{enumerate}
\item The curve $\mathcal{Z}_{h,r}$ has no self-intersections and for any $t\in\mathbb{S}^1$, the limit of the halving chords with direction $u'(t)$ is a chord of $\mathcal{Z}_{h,r}$ that only meets $\mathcal{Z}_{h,r}$ at its endpoints.
\item The curve $\mathcal{Z}_{h,r}$ has Hausdorff dimension equal to the Hausdorff dimension $w$ of the graph of the Weierstrass function $h'$, \ie equal to $w=3-\frac{\ln(\alpha)}{\ln(\beta)}<2$.
\item The curve $\mathcal{Z}_{h,r}$ defines a Lebesgue measurable set with boundary $\mathcal{Z}_{h,r}$ whose Lebesgue measure (area) is
$$A(\mathcal{Z}_{h,r})=\pi\,r^2+\frac{\pi}{2}\,\mathopen{}\left(\frac{1}{\alpha^2-1}-\frac{\beta^2}{\alpha^2-\beta^2}\right)\mathclose{}.$$
\item The curve $\mathcal{Z}_{h,r}$ inherits a Zindler property: its area is split into halves by chords of constant length $2r$.
\end{enumerate}
More specifically, it is enough to take
$r>\frac{1}{\alpha-1}+\frac{2\,\beta}{\alpha-\beta}.$
\end{theorem}

\begin{proof}
\textbf{1.} Define $c_{\alpha,\beta}\deq\frac{1}{\alpha-1}+\frac{\beta}{\alpha-\beta}$. The fractal projective hedgehog $\mathcal{H}_h$ is included in the closed ball centered at $(0,0)$ and with radius $c_{\alpha,\beta}$. Indeed, the distance from a point of $\mathcal{H}_h$ to $(0,0)$ is
\begin{align*}
\bigl\|x_h(t)\bigr\|=\sqrt{h(t)^2+h'(t)^2}&\leq\bigl|h(t)\bigr|+\bigl|h'(t)\bigr|=\mathopen{}\left|\sum_{k=1}^{\infty}\frac{\sin\bigl(\beta^k\,t\bigr)}{\alpha^k}\right|\mathclose{}+\mathopen{}\left|\sum_{k=1}^{\infty}\frac{\beta^k\cos\bigl(\beta^k\,t\bigr)}{\alpha^k}\right|\\
&\leq \sum_{k=1}^{\infty}\mathopen{}\left(\frac{1}{\alpha}\right)^{k}\mathclose{}+\sum_{k=1}^{\infty}\mathopen{}\left(\frac{\beta}{\alpha}\right)^k\mathclose{}=\frac{1}{\alpha-1}+\frac{\beta}{\alpha-\beta}=c_{\alpha,\beta}.
\end{align*}
Moreover, if $r>\frac{1}{\alpha-1}+\frac{2\,\beta}{\alpha-\beta}$, we have
$$r+h'(t)=r+\sum_{k=1}^{\infty}\frac{\beta^k\cos\bigl(\beta^k\,t\bigr)}{\alpha^k}\geq r-\sum_{k=1}^{\infty}\mathopen{}\left(\frac{\beta}{\alpha}\right)^k\mathclose{}=r-\frac{\beta}{\alpha-\beta}>\frac{1}{\alpha-1}+\frac{\beta}{\alpha-\beta}=c_{\alpha,\beta}.$$
Therefore,
$$\bigl\|z_{h,r}(t)\bigr\|=\sqrt{h(t)^2+\bigl(r+h'(t)\bigr)^2}\geq\bigl|r+h'(t)\bigr|>c_{\alpha,\beta},$$
which implies that the fractal curve $\mathcal{Z}_{h,r}$ is outside the closed ball centered at $(0,0)$ and with radius $c_{\alpha,\beta}$.

Now, as a consequence of Theorem~1 of \cite{MartinezMaure00}, we have that for a hedgehog with a $\mathcal{C}^2$-support function, the number of its cooriented support lines passing through a point exterior to the hedgehog is exactly two. Recall that the exterior of a hedgehog $\mathcal{H}_h$ is the set of points $x\in\mathbb{R}^2\setminus\mathcal{H}_h$ such that the winding number of $\mathcal{H}_h=x_h(\mathbb{S}^1)$ around $x$ is zero. We can apply this result to $\mathcal{H}_{h_n}$, because $h_n$ is $\mathcal{C}^2$. Since this hedgehog is projective, we will have that for each point exterior to $\mathcal{H}_{h_n}$, the number of uncooriented support lines is exactly one. We claim that the same property will hold for $\mathcal{H}_h$. Indeed, assume that $\mathcal{H}_h$ admits two different uncooriented support lines, say $L_1$ and $L_2$, that intersect transversely in a point exterior to $\mathcal{H}_h$. These lines $L_1$ and $L_2$ have equations of the form $\langle x, u_1\rangle=h(u_1)$ and $\langle x,u_2\rangle=h(u_2)$ with $u_1$ distinct from $u_2$ and $-u_2$ on $\mathbb{S}^1$. Since $h_n(u_1)$ tends to $h(u_1)$ and $h_n(u_2)$ tends to $h(u_2)$ as $n$ tends to infinity, then the support lines of $\mathcal{H}_{h_n}$, say $L_1(n)$ and $L_2(n)$, with equations $\langle x, u_1\rangle=h_n(u_1)$ and $\langle x,u_2\rangle=h_n(u_2)$, which are different (note they are parallel to $L_1$ and $L_2$, respectively), must intersect in a point exterior to $\mathcal{H}_{h_n}$ when $n$ is large enough. This leads to a contradiction.

Thus, by any point of the fractal curve $\mathcal{Z}_{h,r}$ it passes exactly an uncooriented support line of the fractal projective hedgehog $\mathcal{H}_h$. Therefore, $\mathcal{Z}_{h,r}$ has no self-intersections and any limit halving chord (which is part of such an uncooriented support line) meets the fractal curve $\mathcal{Z}_{h,r}$ only at its endpoints.

\textbf{2.} From the previous item, if $r$ is sufficiently large, then the curve $\mathcal{Z}_{h,r}$ is a Jordan curve (that is, a plane simple closed curve). Therefore, from the Jordan curve theorem, it divides the plane into exactly two regions: a bounded interior and an unbounded exterior. Besides, the Jordan--Sch\"onflies theorem asserts that if $C$ is a Jordan curve, then there exists a homeomorphism $f:\mathbb{R}^2\to\mathbb{R}^2$ such that $f(C)$ is the unit circle in the plane. Therefore, the set $K_{h,r}$ defined by the interior of the curve $\mathcal{Z}_{h,r}$ plus its boundary is homeomorphic to the unit disk and thus compact.

Given $r>0$, let $\tilde{X}:\mathbb{S}^1\times\mathbb{R}\to \mathbb{R}^2$ be the $\mathcal{C}^1$ map defined by
$$\tilde{X}(t,y)=h(t)\,u(t)+(y+r)\,u'(t).$$
The map $\tilde{X}$ is $\mathcal{C}^1$ because its partial derivatives,
\begin{equation*}
\frac{\partial \tilde{X}}{\partial t}(t,y) =\bigl(h'(t)-y-r\bigr)\,u(t) +h(t)\,u'(t) \quad\text{and}\quad \frac{\partial\tilde{X}}{\partial y}(t,y) =u'(t)
\end{equation*}
are continuous on $\mathbb{S}^1\times\mathbb{R}$.

Now, given $r>0$, the function $h'+r$ is continuous on $\mathbb{S}^1$, so it is bounded on $\mathbb{S}^1$. In fact, we have
$$r-\frac{\beta}{\alpha-\beta}\leq r+h'(t)\leq r+\frac{\beta}{\alpha-\beta},$$
where, from what we have seen above, if we take $r>\frac{1}{\alpha-1}+\frac{2\beta}{\alpha-\beta}$, the lower bound is greater than $c_{\alpha,\beta}$.

Given $r>\frac{1}{\alpha-1}+\frac{2\beta}{\alpha-\beta}$, let $X:D\to X(D)\subset\mathbb{R}^2$ be the restriction of $\tilde{X}$ to the compact set $D=\mathbb{S}^1\times\mathopen{}\left[-\frac{\beta}{\alpha-\beta},\frac{\beta}{\alpha-\beta}\right]\mathclose{}$. For the chosen value of $r$, the map $X$ is a $\mathcal{C}^1$-diffeomorphism, where the differential of $X$ at the boundary points of $D$ is defined as the differential of $\tilde{X}$ at these points. Indeed, the map is bijective by the fact that the image $X(D)$ is composed of segments of support lines to $\mathcal{H}_h$ which do not intersect outside the closed ball centered at $(0,0)$ and with radius $c_{\alpha,\beta}$. 

With respect to the standard basis $\mathcal{B}$ in the source and the orthonormal basis $\mathcal{B}(t)\deq \bigl(u(t), u'(t)\bigr)$ in the target, the matrix of $dX_{(t,y)}$ is
\begin{equation*}
\operatorname{Mat}\bigl( dX_{(t,y)};\mathcal{B},\mathcal{B}(t) \bigr) =
{\setlength{\arraycolsep}{0.5em}
\begin{pmatrix}
h'(t)-y-r & 0\\[0.5em]
h(t) & 1
\end{pmatrix},}
\end{equation*}
and its determinant is $h'(t)-y-r$. Since $y\geq -\frac{\beta}{\alpha-\beta}$ and $h'(t)\leq\frac{\beta}{\alpha-\beta}$, we have that
$$h'(t)-y-r\leq h'(t)-r+\frac{\beta}{\alpha-\beta}\leq\frac{2\beta}{\alpha-\beta}-r<\frac{1}{1-\alpha}<0.$$
Therefore, the matrix of $dX_{(t,y)}$ is invertible. Since $D$ is compact and $X:D\to X(D)$ is continuous and bijective, $X^{-1}$ is continuous. With this, and by the inverse function theorem, we deduce that the inverse $X^{-1}$ is also $\mathcal{C}^1$. As a consequence, $X$ is a $\mathcal{C}^1$-diffeomorphism. In fact, if $U(D)$ is a sufficiently small open neighborhood of $D$, the restriction (still denoted by $X$) of $\tilde{X}$ to $U(D)$ is a $\mathcal{C}^1$-diffeomorphism from $U(D)$ onto $X\bigl(U(D)\bigr)$ and $X\bigl(U(D)\bigr)$ is disjoint from the closed ball centered at $(0,0)$ and with radius $c_{\alpha,\beta}$.

Thus, since $X$ is a $\mathcal{C}^1$-diffeomorphism, we have that $X$ is locally bi-Lipschitz. Let
$$\Gamma=\bigl\{(t,h'(t)) :\ t\in\mathbb{S}^1\bigr\},$$
which is the graph of a Weierstrass function $h'$ on $\mathbb{S}^1\cong\mathbb{R}/2\pi\mathbb{Z}$, so it has Hausdorff dimension $w=3-\frac{\ln(\alpha)}{\ln(\beta)}$.

Now, since the Hausdorff dimension of a set is preserved by locally bi-Lipschitz functions (see \eg \cite{FalconerBook, Aberti06}), we have that $\mathcal{Z}_{h,r}=X(\Gamma)$ has Hausdorff dimension $w$.

\textbf{3.} If $r$ is sufficiently large, we have seen before that the set $K_{h,r}$ defined by the interior of the curve $\mathcal{Z}_{h,r}$ plus its boundary is compact. As a consequence, $K_{h,r}$ is a Borel set and therefore a Lebesgue measurable set. Thus, the area of $\mathcal{Z}_{h,r}$ is well-defined by the Lebesgue measure of $K_{h,r}$. Let $K_{h_n,r}$ be the compact set defined by the interior of the curve $\mathcal{Z}_{h_n,r}$ plus its boundary. We have that the sequence of compact sets $K_{h_n,r}$ converges to a compact set $K_{h,r}$ in the Hausdorff metric and the sequence of the corresponding boundaries $\mathcal{Z}_{h_n,r}$ converges to the boundary $\mathcal{Z}_{h,r}$ in the Hausdorff metric. In addition, by the previous item, the Hausdorff dimension of $\mathcal{Z}_{h,r}$ is strictly less than 2. Thus, we can apply Theorem~\ref{th:continuity-lebesgue-measure} and deduce that the sequence of Lebesgue measures of $K_{h_n,r}$ converges to the Lebesgue measure of $K_{h,r}$. Using now Proposition~\ref{prop:area-seq-zindler-curves}, we have
$$A(\mathcal{Z}_{h,r})=\lim_{n\to\infty}\mathopen{}\left[\pi\,r^2+\frac{\pi}{2}\,\mathopen{}\left(\frac{\frac{1}{\alpha^{2n}}-1}{1-\alpha^2}-\frac{\beta^2}{\alpha^2-\beta^2}\,\mathopen{}\left(1-\mathopen{}\left(\frac{\beta}{\alpha}\right)^{2n}\mathclose{}\right)\mathclose{}\right)\mathclose{}\right]\mathclose{}.$$
Since $\alpha>1$ and $\beta/\alpha<1$, we deduce the expression of the statement:
$$A(\mathcal{Z}_{h,r})=\pi\,r^2+\frac{\pi}{2}\,\mathopen{}\left(\frac{1}{\alpha^2-1}-\frac{\beta^2}{\alpha^2-\beta^2}\right)\mathclose{}.$$

\textbf{4.} First of all, notice that if $r>\frac{1}{\alpha-1}+\frac{2\beta}{\alpha-\beta}$, then each $\mathcal{Z}_{h_n,r}$ has no self-intersections and any halving chord meets the curve only at its endpoints. This is because we have similar bounds to those given above:
$$\bigl\|x_{h_n}(t)\bigr\|\leq c_{\alpha,\beta},\qquad\text{and}\qquad \bigl\|z_{h_n,r}(t)\bigr\|>c_{\alpha,\beta},$$ which imply that each $\mathcal{H}_{h_n}$ (resp. $\mathcal{Z}_{h_n,r}$) is also inside (resp. outside) the closed ball centered at $(0,0)$ and with radius $c_{\alpha,\beta}$.

By construction, the chords $C(t)=\mathopen{}\left[z_{h,r}(t),\,z_{h,r}(t+\pi)\right]\mathclose{}$ are of constant length $2r$. Indeed,
$$\bigl\|z_{h,r}(t+\pi)-z_{h,r}(t)\bigr\|=2\,r$$
because $h(t)+h(t+\pi)=0$. In addition, the boundary of $K_{h,r}$ and each chord $C(t)$ splits the set $K_{h,r}$ into two sets $K_{h,r}^1(t)$ and $K_{h,r}^2(t)$ with a well-defined area if $r$ is sufficiently large (these chords only meet $\mathcal{Z}_{h,r}$ at the endpoints). The same can be done for each curve $\mathcal{Z}_{h_n,r}$ and the corresponding positions of the constant length chords $C_n(t)=\mathopen{}\left[z_{h_n,r}(t),\,z_{h_n,r}(t+\pi)\right]\mathclose{}$, where the sets $K_{h_n,r}$ are split into two sets $K_{h_n,r}^1(t)$ and $K_{h_n,r}^2(t)$. The chords $C_n(t)$ tend to $C(t)$ when $n$ tends to infinity. Since all these sets are compact and the boundary of all $K_{h,r}^1(t)$ and $K_{h,r}^2(t)$ also have Hausdorff dimension strictly less than $2$, we can apply again Theorem~\ref{th:continuity-lebesgue-measure} to ensure that
$$\lim_{n\to+\infty}A\bigl(K_{h_n,r}^1(t)\bigr)=A\bigl(K_{h,r}^1(t)\bigr)\qquad\text{and}\qquad\lim_{n\to+\infty}A\bigl(K_{h_n,r}^2(t)\bigr)=A\bigl(K_{h,r}^2(t)\bigr),$$
for all $t\in\mathopen{}\left[0,\pi\right]\mathclose{}$, where $A(K)$ stands for the area given by the positively oriented boundary curve of a compact set $K$. Therefore, since for each $n\in\mathbb{N}^{\ast}$, $\mathcal{Z}_{h_n,r}$ is a Zindler curve with halving chords $C_n(t)$, we have $A\bigl(K_{h_n,r}^1(t)\bigr)=A\bigl(K_{h_n,r}^2(t)\bigr)$ for all $t\in\mathbb{S}^1$. Thus, we deduce $A\bigl(K_{h,r}^1(t)\bigr)=A\bigl(K_{h,r}^2(t)\bigr)$ for all $t\in\mathbb{S}^1$, which implies that $\mathcal{Z}_{h,r}$ is a Zindler curve with $C(t)=\mathopen{}\left[z_{h,r}(t),\,z_{h,r}(t+\pi)\right]\mathclose{}$ being halving chords for its area.
\end{proof}


Since we are on a fractal framework, it is not totally immediate that our fractal Zindler curve $\mathcal{Z}_{h,r}$ with area-halving chords of constant length constitutes indeed the boundary of a fractal body of density $\frac{1}{2}$ that floats in water in equilibrium in every orientation. It is left to prove that, for all orientations, the line that joins the centers of mass of the immersed and emerged parts is perpendicular to the waterline \cite{Ryabogin22annals}.

\begin{theorem}\label{th:fractal-zindler-floating-bodies}
Let $\alpha$, $\beta$, $r$ and $\mathcal{Z}_{h,r}$ be under the assumptions of Theorem~\ref{th:main-theorem-fractal-zindler}. The curve $\mathcal{Z}_{h,r}$ is the boundary of a fractal body of uniform density $\frac{1}{2}$ that floats in water in equilibrium in every orientation.
\end{theorem}

\begin{proof}
The center of mass (centroid) of a Lebesgue measurable compact body $K$ in $\mathbb{R}^2$ is defined by
$$C(K)=\frac{1}{A(K)}\int_K x\,dA,$$
where $A(K)$ is the area of $K$. Notice that this definition is valid for the compact body $K_{h,r}$ defined by our fractal curve $\mathcal{Z}_{h,r}$ as its boundary.

We have that each $\mathcal{Z}_{h_n,r}$ is a smooth Zindler curve (with its boundary defining a compact set $K_{h_n,r}$). Therefore, for each waterline (halving chord) with direction $u'(t)$ in $\mathcal{Z}_{h_n,r}$, the vector $v_{n,r}(t)$ joining the centers of mass of the immersed and emerged parts, say $c_{h_n,r}^1(t)$ and $c_{h_n,r}^2(t)$, is parallel to $u(t)$. 

It is enough to prove that the sequences of centroids $\bigl(c_{h_n,r}^1(t)\bigr)_{n=1}^{\infty}$ and $\bigl(c_{h_n,r}^2(t)\bigr)_{n=1}^{\infty}$ converge to the centroids $c_{h,r}^1(t)$ and $c_{h,r}^2(t)$, respectively, of the immersed and emerged parts of $K_{h,r}$ for each waterline with direction $u'(t)$. Indeed, in such a case, the vector
$$v(t)\deq c_{h,r}^2(t)-c_{h,r}^1(t)=\lim_{n\to\infty}\bigl(c_{h_n,r}^2(t)-c_{h_n,r}^1(t)\bigr)=\lim_{n\to\infty} v_{n,r}(t)$$
is parallel to $u(t)$ for each $t\in\mathbb{S}^1$.

Let $(K_n)_{n=1}^{\infty}$ be a sequence of Lebesgue measurable compact sets converging to a compact set $K$ under the assumptions of Theorem~\ref{th:continuity-lebesgue-measure}. We have
$$\bigl\|\mathds{1}_{K_{n}}-\mathds{1}_{K}\bigr\|_{L^1(\mathbb{R}^2)}=\bigl\|\mathds{1}_{K_{n}\triangle K}\bigr\|_{L^1(\mathbb{R}^2)}=\int_{\mathbb{R}^2} |\mathds{1}_{K_{n}\triangle K}|=\int_{K_{n}\triangle K}1=\lambda(K_{n}\triangle K),$$
where $\mathds{1}_A$ stands for the characteristic function of a set $A\subseteq\mathbb{R}^2$. Since $\lambda(K_{n}\triangle K)$ tends to zero as $n$ goes to infinity, then the sequence $(\mathds{1}_{K_{n}})_{n=1}^{\infty}$ converges pointwise to $\mathds{1}_{K}$ in $L^1(\mathbb{R}^2)$.

Moreover, since all $K_{n}$ are bounded sets, there exists $R>0$ such that all $K_{n}$ are contained in a closed ball $B_R$ centered at the origin and with radius $R$. Hence, given $(x_1,x_2)\in K_{n}$, we have
$$\bigl| x_i\,\mathds{1}_{K_{n}}(x_1,x_2)\bigr|=\bigl|x_i\,\mathds{1}_{K_{n}}(x_1,x_2)\,\mathds{1}_{B_R}(x_1,x_2)\bigr|\leq R\cdot\mathds{1}_{B_R}(x_1,x_2)\in L^1(\mathbb{R}^2),$$
for $i=1,2$.

Thus, we can apply the dominated convergence theorem to each coordinate function of the integrand to ensure that
$$\lim_{n\to\infty}\frac{1}{A(K_{n})}\int_{K_{n}}x\,dA=\lim_{n\to\infty}\frac{1}{A(K_{n})}\int_{\mathbb{R}^2} x\,\mathds{1}_{K_{n}}(x)\,dA=\frac{1}{A(K)}\int_{\mathbb{R}^2} x\,\mathds{1}_{K}(x)\,dA=C(K),$$
where we have used that the sequence of areas $\bigl(A(K_{n})\bigr)_{n=1}^{\infty}$ converges to $A(K)$.

The argument above can be applied to both the sequences $\bigl(K_{h_n,r}^1(t)\bigr)_{n=1}^{\infty}$ and $\bigl(K_{h_n,r}^2(t)\bigr)_{n=1}^{\infty}$ of the immersed and emerged parts, respectively, for each waterline with direction $u'(t)$. 
\end{proof}

\section{Hausdorff measure preservation property of fractal Zindler curves}\label{sec:hauss-measure-zindler}
For a smooth Zindler curve, the chords of constant length bisecting the area are the same as those bisecting the curve's length. In the fractal setting, however, fractal Zindler curves have infinite length, so a direct extension of this definition requires replacing the length of the curve with its $w$-dimensional Hausdorff measure, where $w$ denotes the Hausdorff dimension of the curve. This naturally raises the question: do all halving chords of $\mathcal{Z}_{h,r}$ divide the $w$-dimensional Hausdorff measure of $\mathcal{Z}_{h,r}$ into halves?  As we will prove in this section, the answer is affirmative provided that $r$ is sufficiently large (Corollary~\ref{cor:equality-hauss-measures}).

Henceforth, we assume that $r$ is sufficiently large so that the 4 assertions of Theorem~\ref{th:main-theorem-fractal-zindler} hold. In order to state the main result, Theorem~\ref{th:hauss-measure-preservation}, we introduce the half-period map $S:\mathcal{Z}_{h,r}\to \mathcal{Z}_{h,r}$ defined by $S\bigl( z_{h,r}(t) \bigr) =z_{h,r}\bigl( t+\pi \bigr) $ for all $t\in 
\mathbb{S}^{1}$. Straightforward calculations give:
\begin{enumerate}[(i)]
\item $S^2\deq S\circ S=\operatorname{Id}_{\mathcal{Z}_{h,r}}$ and hence $S^{-1}=S$.
\item $z_{h,r}(t+\pi)=h(t)\,u(t)+\bigl(h'(t)-r)\,u'(t)$ and
hence $z_{h,r}(t+\pi)=z_{h,r}(t)-2\,r\,u'(t)$ for all $t\in \mathbb{S}^{1}$.
\end{enumerate}
Consequently, $S:\mathcal{Z}_{h,r}\to \mathcal{Z}_{h,r}$
is a homeomorphism and $S\bigl(z_{h,r}(t)\bigr) =z_{h,r}(t)-2\,r\,u'(t)$ for all $t\in \mathbb{S}^{1}$. The main result is the following.

\begin{theorem}\label{th:hauss-measure-preservation}
Let $\alpha$, $\beta$, $r$ and $\mathcal{Z}_{h,r}$ be under the assumptions of Theorem~\ref{th:main-theorem-fractal-zindler}. For every $d\geq 0$ and every
Borel set $A$ in $\mathcal{Z}_{h,r}$, the $d$-dimensional Hausdorff measure $H^{d}$ is such that 
\begin{equation*}
H^{d}\bigl(S(A)\bigr)=H^{d}(A).
\end{equation*}%
In particular, if $w$ denotes the Hausdorff dimension of $\mathcal{Z}_{h,r}$, then 
\begin{equation*}
H^{w}\bigl(S(A)\bigr)=H^{w}(A)
\end{equation*}%
for every Borel set $A$ in $\mathcal{Z}_{h,r}$.
\end{theorem}

As a direct consequence of Theorem~\ref{th:hauss-measure-preservation}, we conclude the following halving property for our fractal curve $\mathcal{Z}_{h,r}$.

\begin{corollary}\label{cor:equality-hauss-measures}
Let $\alpha$, $\beta$, $r$ and $\mathcal{Z}_{h,r}$ be under the assumptions of Theorem~\ref{th:main-theorem-fractal-zindler} and let $w=3-\frac{\ln(\alpha)}{\ln(\beta)}$, $1<w<2$, be the Hausdorff dimension of $\mathcal{Z}_{h,r}$. The halving chords of $\mathcal{Z}_{h,r}$ divide the $w$-dimensional Hausdorff measure of $\mathcal{Z}_{h,r}$ into halves.
\end{corollary}


\noindent\textit{Proof of Theorem~\ref{th:hauss-measure-preservation}.}
Denote by $\Gamma$ the graph of $h'$ in the flat
cylinder $C=\mathbb{S}^{1}\times \mathbb{R}$, that is
\begin{equation*}
\Gamma \deq\left\{ p_{t}:=\left( t,h'(t) \right)
\ :\ \  t\in \mathbb{S}^{1}\right\} 
\text{,}
\end{equation*}
and define $X:C\to \mathbb{R}^{2}$ by $X(t,y)
\deq h(t)\,u(t)+(y+r)\,u'(t)$, so that $X(\Gamma)=\mathcal{Z}_{h,r}$ since $X(p_{t}) =z_{h,r}(t) $ for all $t\in \mathbb{S}^{1}$.

Now, define $R:C\to C$ by $R(t,y) =(t+\pi,-y)$, so that
\begin{equation*}
R(p_{t}) =\left( t+\pi ,-h'(t) \right)
=\left( t+\pi ,h'(t+\pi) \right) =p_{t+\pi}
\end{equation*}
and hence $(X\circ R)(p_{t}) =z_{h,r}(t+\pi)=S\bigl(z_{h,r}(t)\bigr)$ for all $t\in \mathbb{S}^{1}$. We thus have the following commutative diagram%
\begin{equation*}
\begin{array}{rcl}
\Gamma & \overset{X}{\rightarrow } & \mathcal{Z}_{h,r} \\ 
R\downarrow &  & \downarrow S \\ 
\Gamma & \overset{X}{\rightarrow } & \mathcal{Z}_{h,r}%
\end{array}%
\end{equation*}
so that, if $\mathcal{Z}_{h,r}$ is regarded as parameterized by $X:\Gamma
\to \mathcal{Z}_{h,r}$, then $R:\Gamma \to \Gamma $ is the map induced by
the half-period map $S:\mathcal{Z}_{h,r}\to \mathcal{Z}_{h,r}$ on $\Gamma$.

\bigskip
\noindent Under our assumptions, from the proof of Theorem~\ref{th:main-theorem-fractal-zindler}, we know the following lemma.

\medskip
\noindent \textit{Lemma 1.} Let $D=\mathbb{S}^1\times\mathopen{}\left[-\frac{\beta}{\alpha-\beta},\frac{\beta}{\alpha-\beta}\right]\mathclose{}$. The map $X:U(D)\to X\bigl(U(D)\bigr)\subset\mathbb{R}^2$ defined by
$$X(t,y)=h(t)\,u(t)+(y+r)\,u'(t)$$
is a $\mathcal{C}^1$-diffeomorphism, where $U(D)$ is an open neighborhood of $D$.

\bigskip
\noindent \textit{Lemma 2}. There exists a $\mathcal{C}^1$-extension $\overline{S}$ of $S$ to $X\bigl(U(D)\bigr)$ such that:
\begin{enumerate}[(1)]
\item $\overline{S}=S$ on $\mathcal{Z}_{h,r}$;
\item For all $t\in\mathbb{S}^1$, $\left( d\overline{S}\right) _{z_{h,r}(t)}:%
\mathbb{R}^{2}\to \mathbb{R}^{2}$ is the orthogonal reflection
accross the line $\mathbb{R}u'(t)$ spanned by $u'(t)$ in $\mathbb{R}^{2}$.
\end{enumerate}

\noindent \textit{Proof of Lemma 2.} From Lemma 1, we can consider the $\mathcal{C}^{1}$ map
\begin{equation*}
\overline{S}=X\circ R\circ X^{-1}:X\bigl(U(D)\bigr)\to (X\circ R)\bigl(U(D)\bigr)\subset \mathbb{R}^{2}.
\end{equation*}

(1) Let $z\in\mathcal{Z}_{h,r}$. Since $\mathcal{Z}_{h,r}$ is a Jordan curve, there exists a unique $s\in\mathbb{S}^1$ such that $z=z_{h,r}(s)=X(p_s)$. Therefore,
$$\overline{S}(z)=(X\circ R)(p_s)=(S\circ X)(p_s)=S(z).$$

(2) For all $t\in \mathbb{S}^{1}$, we have
\begin{equation*}
\bigl( d\overline{S}\bigr) _{z_{h,r}(t)}=dX_{q_{t}}\circ
dR_{p_{t}}\circ \left( dX_{p_{t}}\right)^{-1},
\end{equation*}
where $q_{t}\deq R(p_{t})$, by the chain rule.

We have seen above (proof of Theorem~\ref{th:main-theorem-fractal-zindler}) that
\begin{equation}\label{eq:dXpt-proof-theorem}
dX_{p_{t}}(a,b)=\frac{\partial X}{\partial t}(p_t)\,a+\frac{\partial X}{\partial y}(p_t)\,b =-r\,a\,u(t)+\bigl(h(t)\,a+b\bigr)\,u'(t)
\end{equation}
for all $(a,b)\in \mathbb{R}^{2}$. Since
\begin{equation*}
dR_{p_{t}}(a,b) =\frac{\partial R}{\partial t}(p_{t})\,a+\frac{\partial R}{\partial y}(p_{t})\,b=a\,(1,0) +b\,( 0,-1)=(a,-b),
\end{equation*}
it follows that
\begin{equation*}
\left( dX_{q_{t}}\circ dR_{p_{t}}\right) (a,b) =dX_{q_{t}}(a,-b) =\frac{\partial X}{\partial t}(q_{t})\,a-\frac{%
\partial X}{\partial y}(q_{t})\,b=r\,a\,u(t) +\bigl(h(t)\,a+b)\,u'(t).
\end{equation*}
Therefore, if $dX_{p_{t}}(a,b)$ is written as $\xi\,u(t) +\xi'\,u'(t)$ in the orthonormal
basis $\mathcal{B}(t) =\bigl( u(t), u'(t)\bigr)$, where $(\xi,\xi') \in 
\mathbb{R}^{2}$, then using (\ref{eq:dXpt-proof-theorem}) we get
\begin{equation*}
\bigl( d\overline{S}\bigr) _{z_{h,r}(t)}\bigl(\xi\,u(t) +\xi'\,u'(t)\bigr) =-\xi\,u(t) +\xi'\,u'(t).
\end{equation*}

\noindent In other words, $~\bigl( d\overline{S}\bigr)_{z_{h,r}(t)}:\mathbb{R}^{2}\rightarrow \mathbb{R}^{2}$ is the orthogonal
reflection across the line $\mathbb{R}u'(t)$ spanned
by $u'(t)$ in $\mathbb{R}^{2}$.\qed

\bigskip
We will apply our Lemma 2. Note that the operator norm of $\bigl( d\overline{%
S}\bigr) _{z_{h,r}(t)}:\mathbb{R}^{2}\rightarrow \mathbb{R}^{2}$
is given by%
\begin{equation*}
\left\Vert \bigl( d\overline{S}\bigr)_{z_{h,r}(t)}\right\Vert
_{\textrm{op}}=1,
\end{equation*}
since $\bigl( d\overline{S}\bigr) _{z_{h,r}(t)}:%
\mathbb{R}^{2}\rightarrow \mathbb{R}^{2}$ is an isometry by our Lemma 2.%
\newline

\noindent \textit{Lemma 3}. (\textit{Arbitrary sharp Lipschitz control near }%
$z_{h,r}(t)$). For every $\varepsilon >0$ and $t\in \mathbb{S}^{1}$, there exists $\rho_{t}>0$ such that: 
\begin{equation*}
\bigl\Vert S(x) -S(y) \bigr\Vert \leq (1+\varepsilon)\,\Vert x-y\Vert,\quad\text{for all }x,y\in
\mathcal{Z}_{h,r}\cap B\bigl(z_{h,r}(t); \rho_{t}\bigr).
\end{equation*}

\noindent \textit{Proof of Lemma 3}. As just seen above, $\bigl\Vert \left( d\overline{S}\right)_{z_{h,r}(t)}\bigr\Vert _{\textrm{op}}=1$, and we know that $d\overline{S}$ is continuous at each $z_{h,r}(t)$. Therefore, for each $t\in\mathbb{S}^1$, there exists $\rho_{t}>0$ such that
\begin{equation*}
\left\Vert \bigl(d\overline{S}\bigr)_{x}\right\Vert_{\textrm{op}}\leq (1+\varepsilon), \quad \text{for all }x\in B\bigl(z_{h,r}(t); \rho_{t}\bigr).
\end{equation*}%
Here, we take $\rho_t$ sufficiently small so that $B(z_{h,r}(t); \rho_t)$ is included in $X\bigl(U(D)\bigr)$.

Let $x,y\in B\bigl(z_{h,r}(t); \rho_{t}\bigr) $. By convexity of $%
B\bigl(z_{h,r}(t); \rho_{t}\bigr)$, the line segment $$\left[ x,y%
\right] \deq\left\{ y+\theta\,(x-y) \in \mathbb{R}^{2}\ :\ \theta \in \left[ 0,1\right] \right\} $$ is contained in $B\bigl(z_{h,r}(t); \rho_{t}\bigr)$, and we have:
\begin{equation*}
\bigl\Vert \overline{S}(x) -\overline{S}(y) \bigr\Vert
=\left\Vert \int\nolimits_{0}^{1}\left( \overline{S}\circ \gamma
_{x,y}\right) ^{\prime }\left( \theta \right) d\theta \right\Vert,
\end{equation*}
where $\gamma _{x,y}:\left[ 0,1\right] \to \mathbb{R}%
^{2},\theta \mapsto y+\theta\,(x-y)$, so that%
\begin{align*}
\bigl\Vert \overline{S}(x) -\overline{S}(y) \bigr\Vert & =\left\Vert \int\nolimits_{0}^{1}\bigl(d\overline{S}\bigr)_{\gamma_{x,y}(\theta)}(x-y)\,d\theta \right\Vert \\ 
& \leq \int\nolimits_{0}^{1}\left\Vert \bigl(d\overline{S}\bigr)_{\gamma _{x,y}(\theta)}\right\Vert _{\textrm{op}}\left\Vert x-y\right\Vert d\theta \leq \left( 1+\varepsilon \right) \left\Vert x-y\right\Vert.
\end{align*}

Using (1) of Lemma 2, we finally obtain%
\begin{equation*}
\bigl\Vert S(x) -S(y) \bigr\Vert \leq (
1+\varepsilon)\,\Vert x-y\Vert,\quad \text{for all }x,y\in
\mathcal{Z}_{h,r}\cap B\bigl( z_{h,r}(t); \rho_{t}\bigr).\eqno\qed
\end{equation*}

We are now ready to prove the $d$-dimensional Hausdorff measure preservation
under $S$ stated in Theorem~\ref{th:hauss-measure-preservation}. Fix $d\geq 0$, a Borel set $A$ in $\mathcal{Z}_{h,r}$ and $%
\varepsilon >0$.

By compactness of the Zindler curve $\mathcal{Z}_{h,r}=z_{h,r}\bigl(\mathbb{S}^{1}\bigr)$ in $\mathbb{R}^{2}$, we can extract a minimal finite open subcover $\{
B_{1},\ldots ,B_{N}\} $ of the open cover $\left\{ B\bigl(z_{h,r}(t); \rho_{t}\bigr) \ :\  t\in \mathbb{S}^{1} \right\}$ constructed by Lemma 3:%
\begin{equation*}
\mathcal{Z}_{h,r}\subset \bigcup\limits_{j=1}^{N}B_{j};
\end{equation*}

\noindent we cannot remove any $B_{j}$ from the finite cover without leaving
at least one point of $\mathcal{Z}_{h,r}$ uncovered (minimality of the cover), and $S:\mathcal{Z}_{h,r}\cap B_{j}\to \mathcal{Z}_{h,r}$ is $(1+\varepsilon)$-Lipschitz for every $%
j\in \{1,\ldots ,N\}$. Put $P_{1}\deq \mathcal{Z}_{h,r}\cap B_{1}$, and for every $%
j\in\{ 1,\ldots ,N\}$,
\begin{equation*}
P_{j}\deq \left( \mathcal{Z}_{h,r}\cap B_{j}\right) \cap \left(
\bigcup\limits_{i<j}B_{i}\right) ^{c},
\end{equation*}
where the complement of a subset $E$ in $\mathbb{R}^{2}$ is
denoted by $E^{c}$. Then, $(P_{1},\ldots ,P_{N}) $ is a Borel
partition of $\mathcal{Z}_{h,r}$: $(P_{1},\ldots,P_{N}) $ is a collection of nonempty mutually disjoint Borel subsets of $\mathcal{Z}_{h,r}$ whose union equals $\mathcal{Z}_{h,r}$. Note that by construction, each restriction $S\bigr|_{P_{j}}:P_{j}\to \mathcal{Z}_{h,r}$ is $%
(1+\varepsilon)$-Lipschitz, for $j\in \{ 1,\ldots
,N\}$.

Now, put $A_{j}\deq A\cap P_{j}$ for all $j\in\{ 1,\ldots ,N\} $.
Since $(P_{1},\ldots ,P_{N})$ is a Borel partition of $\mathcal{Z}_{h,r}$, the
sets $A_{j}$ are pairwise disjoint Borel sets such that%
\begin{equation*}
A=\bigcup\limits_{j=1}^{N}A_{j}.
\end{equation*}
Therefore,%
\begin{equation*}
H^{d}(A) =\sum\limits_{j=1}^{N}H^{d}(A_{j}),
\end{equation*}
where the equality is understood in the extended sense. Applying
the fact (Proposition 3.1 from \cite{FalconerBook}) that every $L$-Lipschitz map $f:E\subset \mathbb{R}^{2}\rightarrow 
\mathbb{R}^{2}$ satisfies
\begin{equation*}
H^{d}\bigl(f(E)\bigr) \leq L^{d}\,H^{d}(E)
\end{equation*}
to the $(1+\varepsilon)$-Lipschitz maps $S\bigr|_{P_{j}}:P_{j}\to \mathcal{Z}_{h,r}$, we obtain
\begin{equation*}
H^{d}\bigl( S(A_{j}) \bigr) =H^{d}\left(
S\bigr|_{P_{j}}(A_{j}) \right) \leq (1+\varepsilon)^{d}%
\,H^{d}(A_{j}), \quad \text{for all }j\in \{
1,\ldots ,N\}.
\end{equation*}

Using finite subadditivity of $H^{d}$, it then comes:
\begin{align*}
H^{d}\bigl(S(A)\bigr) & =H^{d}\mathopen{}\left(
S\mathopen{}\left( \bigcup\limits_{j=1}^{N}A_{j}\right)\mathclose{} \right)\mathclose{} =H^{d}\mathopen{}\left( \bigcup\limits_{j=1}^{N}S(A_{j})
\right)\mathclose{}\\ 
& \leq \sum\limits_{j=1}^{N}H^{d}\bigl( S(A_{j})\bigr)  \leq (1+\varepsilon)^{d}\sum\limits_{j=1}^{N}H^{d}(A_{j})=(1+\varepsilon)^{d}\,H^{d}(A).
\end{align*}

\noindent Since $S^{2}=\operatorname{Id}_{\mathcal{Z}_{h,r}}$, we also have $H^{d}(A)
=H^{d}\bigl(S^{2}(A)\bigr) \leq (1+\varepsilon)^{d}\,H^{d}\bigl(S(A)\bigr)$,
and hence:
\begin{equation*}
(1+\varepsilon)^{-d}\,H^{d}(A) \leq 
H^{d}\bigl(S(A)\bigr) \leq (1+\varepsilon)^{d}\,H^{d}(A).
\end{equation*}
If $H^{d}(A)=0$, this implies $H^{d}\bigl(S(A)\bigr)=0$. If $0<H^{d}(A)<+\infty$, letting $\varepsilon \to 0$ gives
\begin{equation*}
H^{d}\bigl(S(A)\bigr) =H^{d}(A).
\end{equation*}
Finally, if $H^{d}(A)=+\infty $, then $%
(1+\varepsilon)^{-d}\,H^{d}(A)\leq 
H^{d}\bigl(S(A)\bigr) $ implies $H^{d}\bigl(S(A)\bigr)=+\infty $. Thus, in every case, we have 
$H^{d}\bigl(S(A)\bigr)=H^{d}(A)$.\qed


\begin{remark}[On the values taken by the $w$-dimensional Hausdorff measure]
Recall that from the proof of Theorem~\ref{th:main-theorem-fractal-zindler}, if $r$ is sufficiently large, we can define a bijective map $X:\Gamma\to\mathcal{Z}_{h,r}$, where $\Gamma=\bigl\{(t,h'(t)) : t\in\mathbb{S}^1\bigr\}$ is the graph of a Weierstrass function. We proved that $X$ is locally bi-Lipschitz. Since $\Gamma$ and $\mathcal{Z}_{h,r}$ are compact, the map $X$ is in fact (globally) bi-Lipschitz. Since bi-Lipschitz maps preserve the Hausdorff measure up to multiplicative constants (Proposition 3.1 from \cite{FalconerBook}), this implies that the Hausdorff measure type (if its value is either positive, zero or infinite) is preserved by bi-Lipschitz maps. As a consequence, we deduce that the $w$-dimensional Hausdorff measure of $\mathcal{Z}_{h,r}$ is of the same type (positive, zero or infinite) as the $w$-dimensional Hausdorff measure of the graph of the associated Weierstrass function on $\mathbb{S}^1\cong \mathbb{R}/2\pi\mathbb{Z}$.

It is known that the $w$-dimensional Hausdorff measure of the graph of a Weierstrass function over a bounded interval is finite (as a consequence of Corollary~11.2 from \cite{FalconerBook}, as the Weierstrass function is Hölder continuous with exponent $2-w$). There exist Weierstrass-type functions with positive $w$-dimensional Hausdorff measure \cite{PU89-fractal}, but the question of whether the $w$-dimensional Hausdorff measure is positive or zero over a bounded interval for the classical Weierstrass functions seems to remain open \cite{RenShen21}.
\end{remark}

\section{Examples of fractal Zindler curves}\label{sec:examples}
Theorems~\ref{th:main-theorem-fractal-zindler} and \ref{th:fractal-zindler-floating-bodies} allows one to easily construct the boundaries of fractal bodies of uniform density $1/2$ floating in water in equilibrium in every orientation (having, thus, exactly half of the body immersed in water for any position). To do so, it is only a matter of choosing suitable parameters $\alpha$, $\beta$ and a sufficiently large $r$.

\begin{figure}[htb]
\begin{minipage}{0.32\textwidth}
\includegraphics[width=\textwidth]{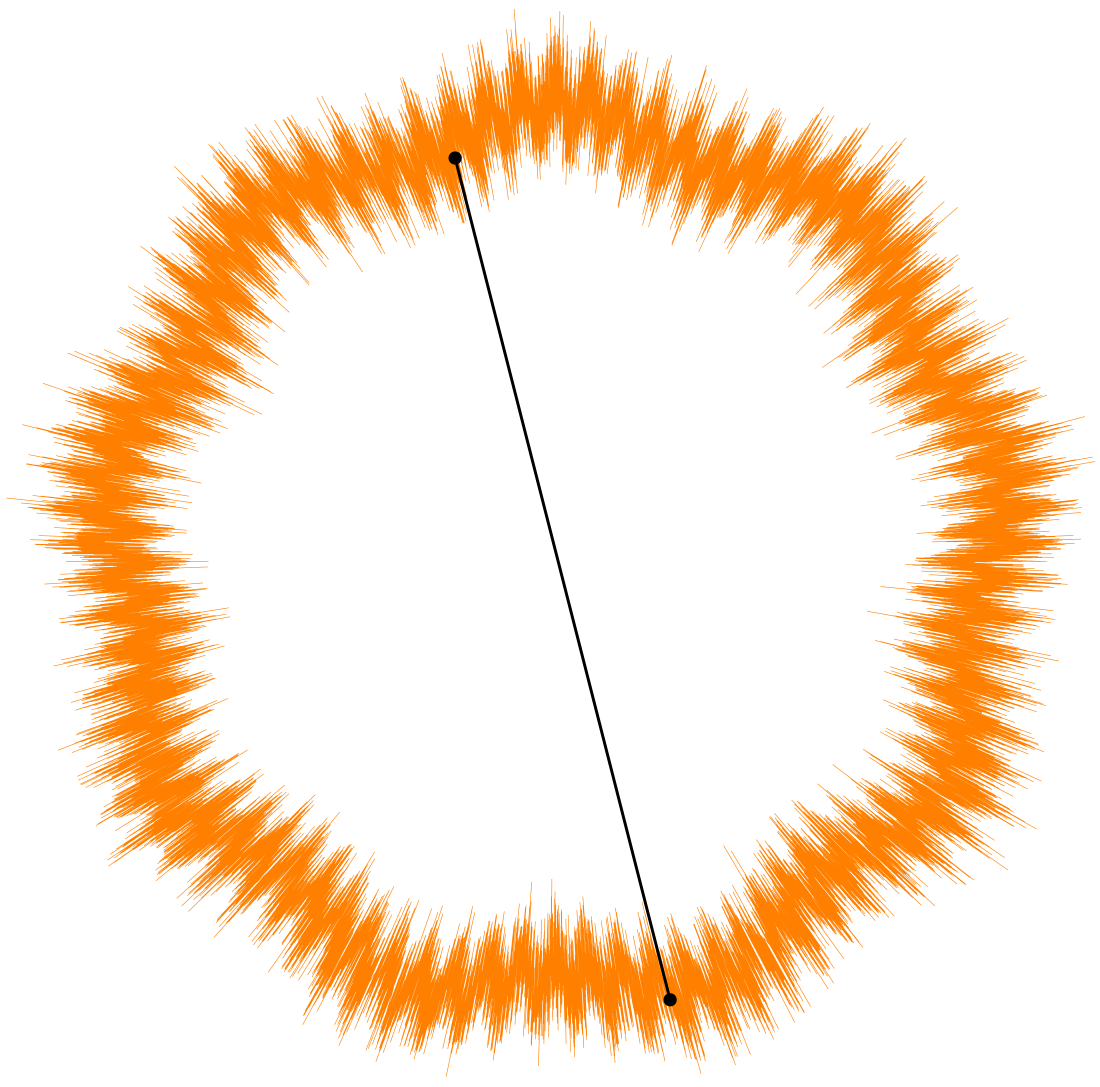}\\[2mm]
\centering {\small \fcolorbox{gray}{white}{$\alpha=10$, $\beta=9$, $r=21$}}
\end{minipage}
\hfill
\begin{minipage}{0.32\textwidth}
\includegraphics[width=\textwidth]{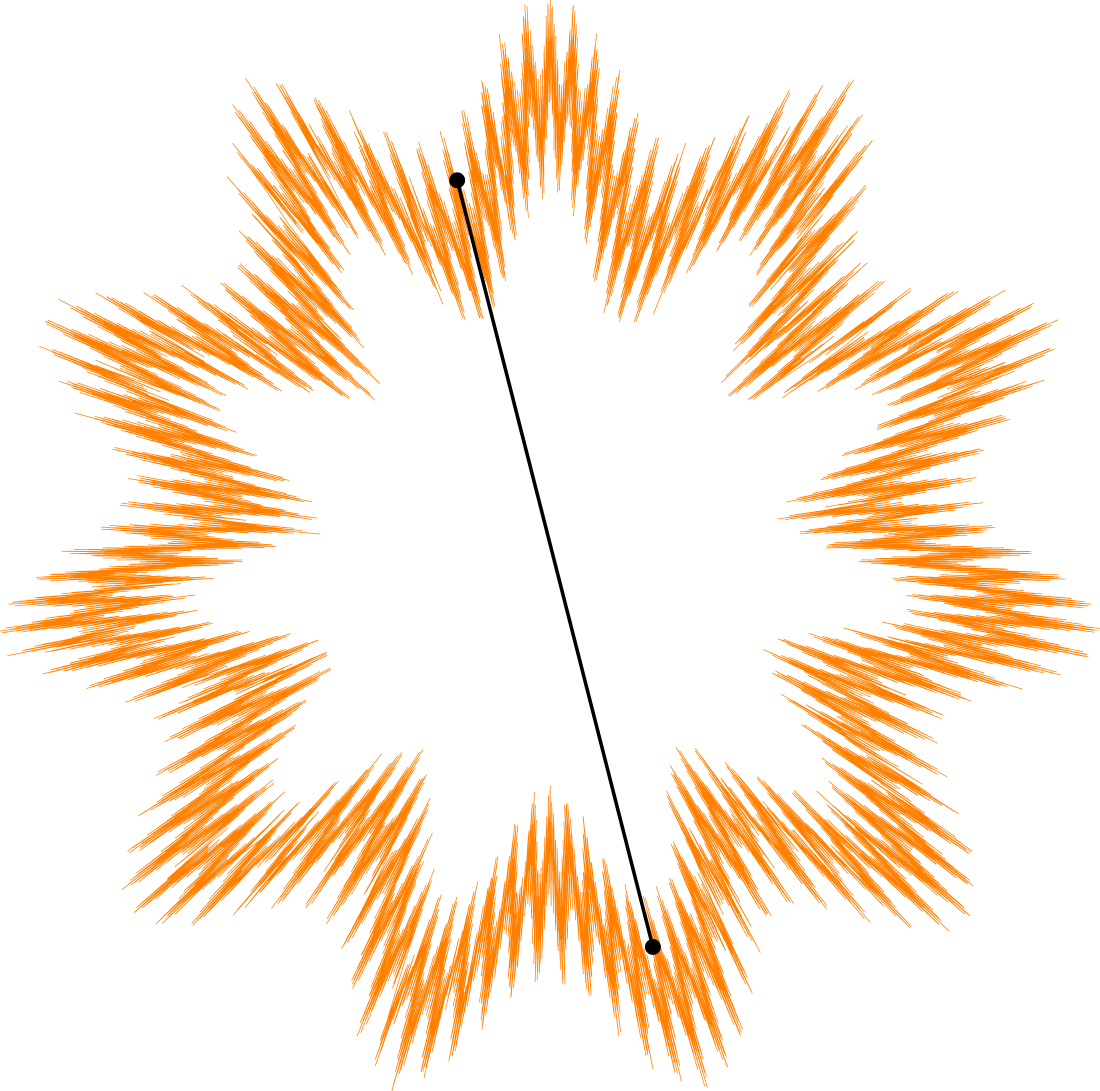}\\[2mm]
\centering {\small \fcolorbox{gray}{white}{$\alpha=17$, $\beta=11$, $r=4$}}
\end{minipage}
\hfill
\begin{minipage}{0.32\textwidth}
\includegraphics[width=\textwidth]{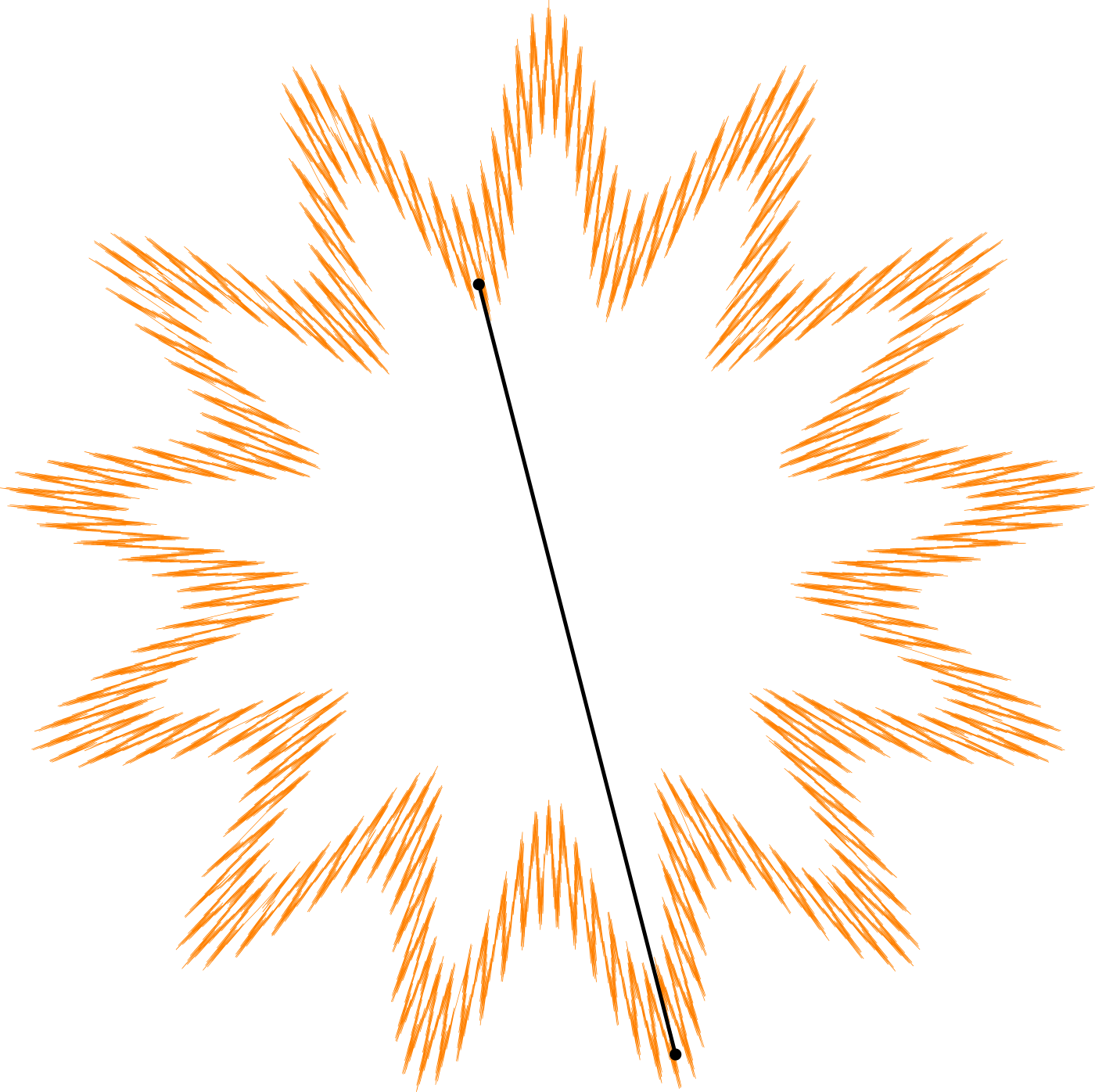}\\[2mm]\centering {\small \fcolorbox{gray}{white}{$\alpha=29$, $\beta=13$, $r=2$}}
\end{minipage}
\caption{Some boundaries $\mathcal{Z}_{h,r}$ of fractal bodies of uniform density $1/2$ floating in equilibrium in every orientation.}
\label{fig:fractal-floating-bodies}
\end{figure}

\figurename~\ref{fig:fractal-floating-bodies} shows some approximations to some fractal Zindler curves that satisfy these conditions.

From a projective hedgehog with support function $h$, it is easy to construct hedgehogs of constant width $2r$ by considering hedgehogs with support functions $h+r$, that is, parallel hedgehogs to $\mathcal{H}_h$ at a signed pointwise distance $r$, measured at each point $x_h(t)$ in the direction of $u(t)$: $x_{h+r}(t)=x_h(t)+r\,u(t)$ for all $t\in\mathbb{S}^1$. Also, we can associate a Zindler curve $\mathcal{Z}_{h,r}$ with halving chords of length $2r$ from the projective hedgehog by offsetting it a distance $r$ in the direction $u'(t)$: $z_{h,r}(t)=x_h(t) + r\,u'(t)$ for all $t\in\mathbb{S}^1$. By construction these two curves are associated through the projective hedgehog and they geometrically correspond to a right-angle rotation of the chords of constant length.

\begin{figure}[htb]
\centering
\begin{overpic}[width=0.7\textwidth]{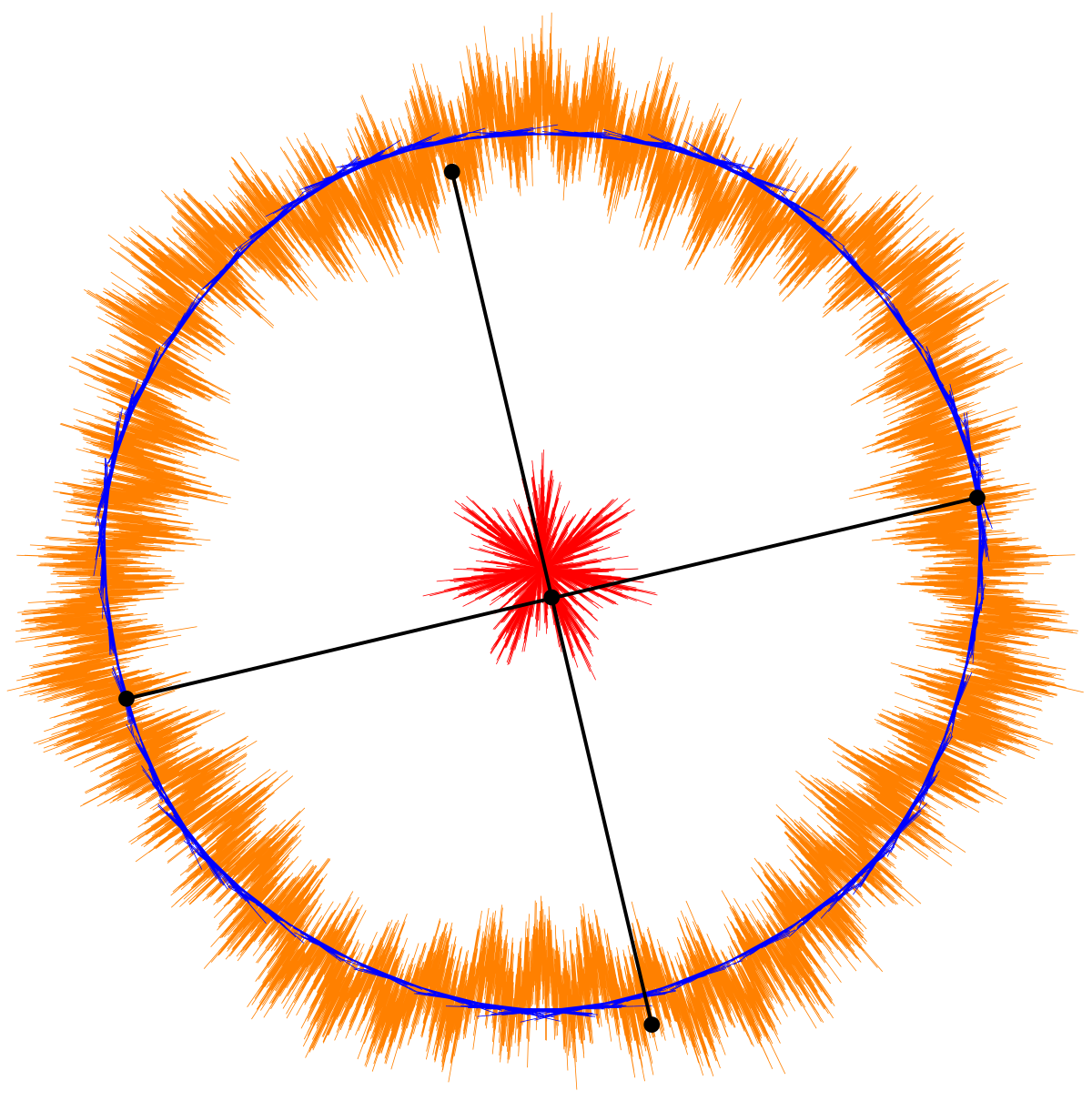}
    \put(72,89){\textcolor{orange}{$\mathcal{Z}_{h,r}$}}
    \put(85,70){\textcolor{blue}{$\mathcal{H}_{h+r}$}}
    \put(56,57){\textcolor{red}{$\mathcal{H}_{h}$}}		    
\end{overpic}
\caption{An approximation (with $n=50$) to a fractal Zindler curve $\mathcal{Z}_{h,r}$ and its corresponding fractal hedgehog of constant width $\mathcal{H}_{h+r}$ via a fractal projective hedgehog $\mathcal{H}_h$, for $\alpha=8$, $\beta=7$ and $r=17$.}
\label{fig:example-fractal-zindler-hcw}
\end{figure}

The same can be repeated in our fractal framework to obtain pairs of associated fractal curves from a fractal projective hedgehog $\mathcal{H}_h$. This way, we can obtain fractal hedgehogs of constant width (that will in fact be the limit of a sequence of smooth hedgehogs of constant width) associated to fractal Zindler curves. In addition, if $r$ is sufficiently large, these fractal Zindler curves are boundaries of figures of uniform density $1/2$ that float in water in equilibrium in every orientation. See an example of this kind in \figurename~\ref{fig:example-fractal-zindler-hcw}.

\appendix
\section{The value of some integrals}\label{sec:appendix}

\begin{lemma*}\label{lem:integrals}
If $n\in\mathbb{N}^{\ast}$, $\beta\in\mathbb{Z}$ and $\alpha\in\mathbb{R}$, $\alpha\neq0$, then
\begin{enumerate}
\item $\displaystyle \int_0^{2\pi}\mathopen{}\left(\sum_{k=1}^n\frac{\sin(\beta^k\,t)}{\alpha^k}\right)^2\mathclose{}dt=\frac{\pi}{\alpha^{2n}}\,\sum_{k=0}^{n-1}\alpha^{2k}.$
\item $\displaystyle \int_0^{2\pi}\sum_{k=1}^n\mathopen{}\left(\frac{\beta}{\alpha}\right)^k\mathclose{}\cos(\beta^k\,t)\,dt=0.$
\item $\displaystyle \int_0^{2\pi}\mathopen{}\left(\sum_{k=1}^n\frac{\sin(\beta^k\,t)}{\alpha^k}\right)\mathclose{}\mathopen{}\left(\sum_{k=1}^n\mathopen{}\left(\frac{\beta^2}{\alpha}\right)^k\mathclose{}\sin(\beta^k\,t)\right)\mathclose{}dt=\frac{\beta^2\,\pi}{\alpha^{2n}}\,\sum_{k=0}^{n-1}\alpha^{2n-2k-2}\,\beta^{2k}.$
\end{enumerate}
\end{lemma*}

\begin{proof}
First, since $\beta\in\mathbb{Z}$, notice that
\begin{equation}\label{eq:integral-zero}
\int_0^{2\pi}\sin(\beta^{n+1}\,t)\sin(\beta^k\,t)\,dt=0,
\end{equation}
for all $n,k\in\mathbb{N}$. The computation of this integral is straightforward by using the identity $\sin(x)\sin(y)=\frac12\,\bigl(\cos(x-y)-\cos(x+y)\bigr)$, that holds for any $x,y\in\mathbb{R}$.

Now, we will prove the first expression of the statement by induction over $n$. For $n=1$, it is easy to compute the integral and since $\beta\in\mathbb{Z}$, the formula is satisfied:
$$\int_0^{2\pi}\mathopen{}\left(\frac{\sin(\beta\,t)}{\alpha}\right)^2\mathclose{}\;dt=\frac{\pi}{\alpha^2}.$$
Suppose that the statement is true for $n$ and let us show it for $n+1$. We have
\begin{align*}
\int_0^{2\pi}&\mathopen{}\left(\sum_{k=1}^{n+1}\frac{\sin(\beta^k\,t)}{\alpha^k}\right)^2\mathclose{}dt=\int_0^{2\pi}\mathopen{}\left(\sum_{k=1}^n\frac{\sin(\beta^k\,t)}{\alpha^k}+\frac{\sin(\beta^{n+1}\,t)}{\alpha^{n+1}}\right)^2\mathclose{}dt\\
&=\int_0^{2\pi}\mathopen{}\left(\mathopen{}\left(\sum_{k=1}^n\frac{\sin(\beta^k\,t)}{\alpha^k}\right)^2\mathclose{}+\frac{\sin^2(\beta^{n+1}\,t)}{\alpha^{2n+2}}+2\,\frac{\sin(\beta^{n+1}\,t)}{\alpha^{n+1}}\,\sum_{k=1}^n\frac{\sin(\beta^k\,t)}{\alpha^k}\right)\mathclose{}dt\\
&=\int_0^{2\pi}\mathopen{}\left(\sum_{k=1}^n\frac{\sin(\beta^k\,t)}{\alpha^k}\right)^2\mathclose{}dt+\int_0^{2\pi}\frac{\sin^2(\beta^{n+1}\,t)}{\alpha^{2n+2}}\,dt+2\,\int_0^{2\pi}\frac{\sin(\beta^{n+1}\,t)}{\alpha^{n+1}}\,\sum_{k=1}^n\frac{\sin(\beta^k\,t)}{\alpha^k}\,dt.
\end{align*}
We use the inductive step for the first integral, the evaluation of the second integral is straightforward and the third integral can be written as
$$\int_0^{2\pi}\frac{\sin(\beta^{n+1}\,t)}{\alpha^{n+1}}\,\sum_{k=1}^n\frac{\sin(\beta^k\,t)}{\alpha^k}\,dt=\sum_{k=1}^n\frac{1}{\alpha^{n+k+1}}\,\int_0^{2\pi}\sin(\beta^{n+1}\,t)\sin(\beta^k\,t)\,dt,$$
which is equal to zero by (\ref{eq:integral-zero}). Therefore, we get
\begin{align*}
\int_0^{2\pi}\mathopen{}\left(\sum_{k=1}^{n+1}\frac{\sin(\beta^k\,t)}{\alpha^k}\right)^2\mathclose{}dt&=\frac{\pi}{\alpha^{2n}}\,\sum_{k=0}^{n-1}\alpha^{2k}+\frac{\pi}{\alpha^{2n+2}}=\frac{\pi}{\alpha^{2n+2}}\,\mathopen{}\left(1+\alpha^2\sum_{k=0}^{n-1}\alpha^{2k}\right)\mathclose{}\\
&=\frac{\pi}{\alpha^{2n+2}}\,\mathopen{}\left(1+\sum_{k=1}^n\alpha^{2k}\right)\mathclose{}=\frac{\pi}{\alpha^{2n+2}}\,\sum_{k=0}^n\alpha^{2k}.
\end{align*}
This completes the proof of the first item. The second item follows immediately by the linearity of the integral and because $\beta\in\mathbb{Z}$. For the third item, we proceed again by induction over $n$. The case $n=1$ is satisfied:
$$\int_0^{2\pi}\frac{\sin(\beta\,t)}{\alpha}\,\mathopen{}\left(\frac{\beta^2}{\alpha}\right)\mathclose{}\sin(\beta\,t)\,dt=\int_0^{2\pi}\frac{\beta^2}{\alpha^2}\sin^2(\beta\,t)\,dt=\frac{\pi\,\beta^2}{\alpha^2}.$$
Now, suppose that the expression is valid for $n$ and we will prove it for $n+1$. We have
\begin{align*}
\int_0^{2\pi}&\mathopen{}\left(\sum_{k=1}^{n+1}\frac{\sin(\beta^k\,t)}{\alpha^k}\right)\mathclose{}\,\mathopen{}\left(\sum_{k=1}^{n+1}\mathopen{}\left(\frac{\beta^2}{\alpha}\right)^k\mathclose{}\sin(\beta^k\,t)\right)\mathclose{}dt\\
&=\int_0^{2\pi}\mathopen{}\left(\sum_{k=1}^n\frac{\sin(\beta^k\,t)}{\alpha^k}+\frac{\sin(\beta^{n+1}\,t)}{\alpha^{n+1}}\right)\mathclose{}\mathopen{}\left(\sum_{k=1}^n\mathopen{}\left(\frac{\beta^2}{\alpha}\right)^k\mathclose{}\sin(\beta^k\,t)+\mathopen{}\left(\frac{\beta^2}{\alpha}\right)^{n+1}\mathclose{}\sin(\beta^{n+1}\,t)\right)\mathclose{}dt\\
&=\int_0^{2\pi}\mathopen{}\left(\sum_{k=1}^n\frac{\sin(\beta^k\,t)}{\alpha^k}\right)\mathclose{}\mathopen{}\left(\sum_{k=1}^n\mathopen{}\left(\frac{\beta^2}{\alpha}\right)^k\mathclose{}\sin(\beta^k\,t)\right)\mathclose{}dt\\
&\qquad+\sum_{k=1}^n\frac{\beta^{2n+2}}{\alpha^{n+k+1}}\int_0^{2\pi}\sin(\beta^k\,t)\sin(\beta^{n+1}\,t)\,dt\\
&\qquad+\sum_{k=1}^n\frac{\beta^{2k}}{\alpha^{n+k+1}}\int_0^{2\pi}\sin(\beta^{n+1}\,t)\sin(\beta^k\,t)\,dt+\frac{\beta^{2n+2}}{\alpha^{2n+2}}\int_0^{2\pi}\sin^2(\beta^{n+1}\,t)\,dt.
\end{align*}
We use the inductive step in the first integral. The two integrals in the middle are zero by (\ref{eq:integral-zero}). The evaluation of the last integral is straightforward. Therefore,
\begin{align*}
\int_0^{2\pi}&\mathopen{}\left(\sum_{k=1}^{n+1}\frac{\sin(\beta^k\,t)}{\alpha^k}\right)\mathclose{}\,\mathopen{}\left(\sum_{k=1}^{n+1}\mathopen{}\left(\frac{\beta^2}{\alpha}\right)^k\mathclose{}\sin(\beta^k\,t)\right)\mathclose{}dt=\frac{\beta^2\pi}{\alpha^{2n}}\sum_{k=0}^{n-1}\alpha^{2n-2k-2}\beta^{2k}+\frac{\beta^{2n+2}\pi}{\alpha^{2n+2}}\\
&=\frac{\beta^2\pi}{\alpha^{2n+2}}\,\mathopen{}\left(\alpha^2\sum_{k=0}^{n-1}\alpha^{2n-2k-2}\beta^{2k}+\beta^{2n}\right)\mathclose{}=\frac{\beta^2\pi}{\alpha^{2n+2}}\,\mathopen{}\left(\sum_{k=0}^{n-1}\alpha^{2n-2k}\beta^{2k}+\beta^{2n}\right)\mathclose{}\\
&=\frac{\beta^2\pi}{\alpha^{2n+2}}\,\mathopen{}\sum_{k=0}^{n}\alpha^{2n-2k}\beta^{2k},
\end{align*}
which completes the proof.
\end{proof}

\section*{Acknowledgements}
The second author has been partially supported by Grant PID2025-172404NB-I00 funded by MCIN/AEI/10.13039/501100011033 and by “ERDF A way of making Europe”, and by Grant UV-INV\_AE\_4225594.

\bibliographystyle{amsplain}
\bibliography{biblio}

\end{document}